\documentclass[twocolumn]{autart}
\usepackage{savesym} \savesymbol{AND}

\usepackage{url}
\usepackage[T1]{fontenc}
\usepackage[latin9]{inputenc}
\usepackage{color}
\usepackage{enumerate,graphicx,epstopdf}
\usepackage{amsmath,amssymb,bm}
\usepackage{cite}
\usepackage{mathtools}
\usepackage{array}
\usepackage{algpseudocode,algorithm,algorithmicx}
\usepackage{color}
\usepackage{lipsum}

\newtheorem{rmk}{Remark}
\newtheorem{dfn}{Definition}
\newtheorem{lmm}{Lemma}
\newtheorem{prp}{Proposition}

\newcommand\reals[0]{\mathbb{R}}

\newcommand\eps[0]{\varepsilon}

\newcommand\dom[1]{\mathrm{dom}~{#1}}
\newcommand{\gph}{\mathrm{gph}~}
\newcommand\dist[2]{\mathrm{dist}({#1},{#2})}
\newcommand\cdiff[0]{\bar{\partial}}
\newcommand\oh[0]{o}
\newcommand\Oh[0]{O}
\newcommand\intrng[2]{\mathbb{Z}_{[#1,#2]}}
\DeclareMathOperator{\diag}{diag}

\begin{document}

\begin{frontmatter}

\title{FBstab: A Stabilized Semismooth Quadratic Programming Algorithm with Applications in Model Predictive Control\thanksref{footnoteinfo}}

\thanks[footnoteinfo]{This research is supported by the National Science Foundation Award Number CMMI 1562209.}

\author[UM]{Dominic Liao-McPherson}\ead{dliaomcp@umich.edu},
\author[UM]{Ilya Kolmanovsky}\ead{ilya@umich.edu}

\address[UM]{Department of Aerospace Engineering, University of Michigan, 1221 Beal Avenue, Ann Arbor, MI 48109 } 


\begin{abstract}
This paper introduces the proximally stabilized Fischer-Burmeister method (FBstab); a new algorithm for convex quadratic programming that synergistically combines the proximal point algorithm with a primal-dual semismooth Newton-type method. FBstab is numerically robust, easy to warmstart, handles degenerate primal-dual solutions, detects infeasibility/unboundedness and requires only that the Hessian matrix be positive semidefinite. We outline the algorithm, provide convergence and convergence rate proofs, report some numerical results from model predictive control benchmarks, and also include experimental results. We show that FBstab is competitive with and often superior to, state of the art methods, has attractive scaling properties, and is especially promising for model predictive control applications.
\end{abstract}

\end{frontmatter}

\section{Introduction}

Model Predictive Control (MPC)\cite{borrelli2017predictive,grune2017nonlinear} is an optimization based control methodology that is gaining in popularity for controlling constrained and/or nonlinear systems. In the case of a linear prediction model, polyhedral constraints and quadratic costs the MPC control law is defined by the solution of a quadratic program (QP). Moreover, QPs are commonly used in methods for nonlinear MPC e.g., in sequential quadratic programming (SQP) \cite{nocedal2006numerical} based methods, such as the real-time iteration scheme \cite{diehl2005real}, for nonlinear MPC, or branch and bound methods for mixed-integer optimization \cite{fletcher1998numerical}. These QPs need to be solved in real-time on embedded systems with limited computing power; warmstarting, where the solver is initialized using a solution from the previous sampling instance, and structure exploitation are often necessary to meet real-time requirements. Moreover, reliability, exception safety, robustness to early termination, and infeasibility detection are important concerns due to the safety critical nature of many MPC controllers.

A number of useful algorithms and packages for solving convex QPs have been developed including: Active Set (AS) methods \cite{ferreau2014qpoases,goldfarb1983numerically,bemporad2016quadratic,cimini2019complexity}, Interior Point (IP) methods \cite{wang2010fast,rao1998application,domahidi2013ecos,mattingley2012cvxgen}, first order (FO) methods \cite{patrinos2014accelerated,patrinos2015dual,stellato2017osqp,goldstein2014fast}, and Dual Newton (DN) methods \cite{patrinos2011global,frasch2015parallel}. AS methods are typically very fast for small to medium problem sizes and can be easily warmstarted using an estimate of the active constraint set. However, they do not scale well, since they have difficulties exploiting sparsity, and are not robust to early termination, i.e., intermediate iterates do not approximate the solution in a meaningful way. IP methods are fast, efficient, robust to early termination and can solve large, structured problems efficiently using advanced linear algebra techniques. However, they are notoriously difficult to warmstart. FO methods can be warmstarted easily, and are attractive from a certification standpoint due to their simplicity and the availability of tight complexity bounds. However, they have slow convergence rates relative to AS, IP, and DN methods; they tend to be most effective on small, strongly convex, or simply constrained problems. Dual Newton methods can exploit structure and be warmstarted but require restrictive assumptions, e.g., that the QP be strongly convex and the linear independence constraint qualification (LICQ) holds. These assumptions reduce their applicability and may cause robustness issues. Finally, primal-dual Newton-type methods, e.g., FBRS\cite{liaomcphersonFBRS}, retain the same warmstarting and structure exploitation as DN methods but relax the strong convexity requirement to the weaker strong second order sufficient condition (SSOSC) but still require the LICQ.

Recently, a hybrid method has been proposed that combines a first order and active set method. QPNNLS\cite{bemporad2018numerically} uses the proximal point algorithm \cite{rockafellar1976monotone} to construct a sequence of regularized QP subproblems whose solutions converge to the solution of the original problem. Each regularized QP is strictly convex and is efficiently solved using a non-negative least squares based active set method \cite{bemporad2016quadratic}. Since proximal-point subproblems are expensive, QPNNLS heavily exploits warmstarting to reduce the cost of solving subsequent subproblems.

In this paper, we propose the proximally stabilized Fischer-Burmeister method (FBstab). FBstab is a hybrid method in the same vein as QPNNLS\cite{bemporad2018numerically} but we employ a primal-dual version of the proximal point algorithm (i.e., the proximal method of multipliers \cite{rockafellar1976augmented}) rather than a primal version, and solve the proximal subproblems with a primal-dual Newton-type method. Using a Newton-type method, rather than an active set method, allows us to solve the proximal subproblems inexactly in addition to warmstarting them, leading to considerable computational savings. In turn, the proximal regularization yields subproblems which automatically satisfy the regularity conditions needed to guarantee robustness and rapid convergence of the Newton-type method. 

The contributions of this paper are as follows:
(i) We describe FBstab, a method which is numerically robust, and can handle problems with degenerate solutions. Moreover, it is easy to warmstart making it compatible with SQP or suboptimal \cite{allan2017inherent} methods and thus ideal for MPC applications.
(ii) We provide detailed convergence and convergence rate proofs under only the assumption that the Hessian is positive semidefinite and a solution exists.
(iii) When a solution does not exist we prove that FBstab can detect and certify infeasibility and unboundedness.
(iv) We illustrate the performance of FBstab through numerical examples and an experimental case study. We also demonstrate that, since the linear systems at the core of FBstab are structured similarly to some IP methods, we can exploit existing specialized linear algebra routines for MPC problems to achieve linear execution time scaling in the receding horizon length.
(v) We have prepared an open source implementation of FBstab which is available at: 
\begin{equation*}
	\text{\url{https://github.com/dliaomcp/fbstab-matlab.git}}
\end{equation*}

FBstab significantly extends FBRS \cite{liaomcphersonFBRS} by removing the SSOSC and LICQ assumptions, improving its numerical robustness, and enabling infeasibility detection, while still being easy to warmstart and retaining the ability to exploit structure. FBstab requires almost no assumptions aside from (non-strict) convexity making it robust and capable of solving any convex QP; only IP methods based on self-dual embedding, e.g., ECOS\cite{domahidi2013ecos}, QPNNLS \cite{bemporad2018numerically} and the alternating direction method of multipliers (ADMM) \cite{goldstein2014fast,stellato2017osqp} are as widely applicable. Finally the use of a primal-dual proximal point method allows FBstab to detect unboundedness, i.e., dual infeasibility, as well as primal infeasibility, unlike QPNNLS which can only detect primal infeasibility. This paper significantly extends the conference version\cite{cdc_fbstab2019}. In particular, it contains detailed convergence, rate-of-convergence, and infeasibility detection proofs, an additional spacecraft relative motion example, an experimental demonstration, and more numerical experiments, all which are not included in \cite{cdc_fbstab2019}.

\section{Problem Setting and Mathematical Preliminaries}
In this paper, we consider convex QPs of the following form,
\begin{subequations}\label{eq:QP}
\begin{align} 
\underset{z}{\mathrm{min.}} \quad \frac12 &z^T H z + f^T z, \\
\mathrm{s.t} ~\quad &Gz = h,\\
\quad &Az \leq b,
\end{align}
\end{subequations}
where $H = H^T \succeq 0 \in \reals^{n\times n}$ is the Hessian matrix, $f \in \reals^n$, $z\in \reals^n$, $G \in \reals^{m \times n}$, $h\in \reals^m$, $A\in\reals^{q \times n}$, and $b \in \reals^q$. We make no assumptions about the problem data aside from the positive semidefiniteness of the Hessian. The Lagrangian for this problem is
\begin{equation*}
  L(z,\lambda,v) = \frac12 z^T H z + f^T z + \lambda^T(Gz-h) + v^T(Az-b),
\end{equation*}
where $\lambda \in \reals^m$ and $v\in \reals^q$ are dual variables, and its dual is
\begin{subequations}\label{eq:dual_QP}
\begin{gather} 
\underset{u,\lambda,v}{\mathrm{min.}} \quad \frac12 u^T H u + b^Tv + h^T\lambda, \\
\mathrm{s.t} \quad Hu + f + G^T \lambda + A^T v = 0,\\
\quad v \geq 0.
\end{gather}
\end{subequations}

We will use $x = (z,\lambda,v) \in \reals^l$ to denote the primal-dual triple. The Karush-Kuhn-Tucker (KKT) conditions for the problem are 
\begin{subequations} \label{eq:KKT}
\begin{gather}
\nabla_z L(z,\lambda,v) = 0,\\
Gz = h,\\
Az - b \leq 0,~v \geq 0,~v^T (Az-b) = 0. \label{eq:cmp_end}
\end{gather}
\end{subequations}
Any vector satisfying \eqref{eq:KKT} is called a critical point. If the feasible set 
\begin{equation}
	\Gamma = \{z~|~Az \leq b,~Gz = h\},
\end{equation}
is nonempty, then the KKT conditions are necessary and sufficient for global optimality \cite{boyd2004convex}. The objective of the FBstab algorithm is to either find a vector $x^*$ which satisfies \eqref{eq:KKT} or to prove that no such vector exists. It does this by applying Newton's method to a semismooth reformulation of \eqref{eq:KKT}. The proximal point algorithm is then wrapped around the semismooth Newton's method; it stabilizes the Newton iterations and detects primal-dual infeasibility. This improves the robustness of the algorithm and allows it to solve problems with degenerate primal-dual solutions, a class of problems which are typically challenging for Newton-type methods.

\subsection{Nonsmooth Analysis} In this section we review some concepts from non-smooth analysis which are used to construct the core Newton-type method. Suppose a function $G:\reals^N \to \reals^M$ is locally Lipschitz on a set $U \subseteq \reals^N$, so that $G$ is differentiable almost everywhere by Rademacher's theorem \cite{rademacher1919partielle}. Clarke's generalized Jacobian \cite{clarke1990optimization} is defined as
\begin{multline}
  \partial G(x) = \mathrm{co}~\{J\in \reals^{M\times N}~| \\~ \exists \{x^k\} \subset D_G : \{x^k\} \rightarrow x,~ \{\nabla G(x_k)\} \rightarrow J\},
\end{multline}
where $D_G$ is the dense set of points where $G$ is differentiable, and co $(\cdot)$ denotes the convex hull. 

We also make use of the C-subdifferential \cite{qi1996c} which is defined as
\begin{equation}
  \bar{\partial} G = \partial G_1 \times \partial G_2 \times ... \times~\partial G_{M},
\end{equation}
where $\partial G_i$ are the generalized Jacobians of the components mappings of $G$. Note that each element of $\partial G_i$ is a row vector. This form of the C-subdifferential is used in \cite{chen1998global,chen2000penalized} and possesses many of the useful properties of the generalized Jacobian but is easier to compute and characterize. 

A function $G:\reals^N \mapsto \reals^M$ is said to be semismooth \cite{qi1993nonsmooth} at $x \in \reals^N$ if $G$ is locally Lipschitz at $x$, directionally differentiable in every direction and the estimate
\begin{equation} \label{eq:genjac_sup}
  \underset{J \in \partial G(x+\xi)}{\text{sup}} ||G(x+\xi) - G(x) - J\xi|| = o(||\xi||),
\end{equation}
holds\footnote{See \cite[A.2]{izmailov2014newton} or \cite[A.2]{nocedal2006numerical} for more details on O notation.}. If $o(||\xi||)$ is replaced with $O(||\xi||^2)$ in \eqref{eq:genjac_sup} then $G$ is said to be strongly semismooth at $x$. The generalized Jacobian and the C-subdifferential can be used to construct Newton-type methods for semismooth systems of nonlinear equations \cite{qi1993nonsmooth,qi1996c}. The following Newton-type method,
\begin{equation}
  x_{k+1} = x_k - V^{-1}G(x_k), \quad V \in \partial G(x_k),
\end{equation}
is locally superlinearly convergent to roots of $G$ which satisfy some regularity properties \cite{qi1993nonsmooth,qi1993convergence}. Similar results are available using the C-subdifferential\cite{qi1996c,chen2000penalized}.

\subsection{Monotone Operators}
In this section we review the proximal point algorithm which is used to stabilize FBstab's core Newton-type method. Recall that a set-valued mapping\footnote{Equivalently a multifunction, operator, relation, point to set mapping, or correspondence.} $T:D \rightrightarrows \reals^N$ is said to be monotone if
\begin{equation}
  \langle x - y, u - v\rangle \geq 0 ,\quad \forall~u \in T(x), v \in T(y), ~ x,y \in D,
\end{equation}
where $D = \dom{T} \subseteq \reals^N$. In addition, if 
\begin{equation}
  \gph{T} = \{(x,u) \in D \times D ~|~ u\in T(x)\},
\end{equation}
is not properly contained in the graph of any other monotone operator then $T$ is said to be maximal \cite{rockafellar1976monotone}. A useful example of a maximal monotone operator is the normal cone. The normal cone mapping of a nonempty closed convex set $C \subseteq \reals^N$ is defined by
\begin{equation*}
 N_C(x) = \begin{cases}
  \{w~|~\langle x - u,w\rangle \geq 0,~\forall u \in C\}, & x\in C,\\
  \emptyset, & x \notin C.
\end{cases}
\end{equation*}

The proximal point algorithm \cite{rockafellar1976monotone} can be used to find zeros of maximal monotone operators. This algorithm  generates a sequence $\{x_k\}$ by the rule
\begin{equation} \label{eq:prox}
  x_{k+1} = P_k(x_k), \quad P_k = (I + \sigma_k^{-1} T)^{-1},
\end{equation}
where $\sigma_k$ is a sequence of positive numbers. The parameter $\sigma_k$ acts as a regularization term; an advantage of the proximal point method over standard regularization techniques is that $\sigma_k$ need not to be driven to zero. Since $T$ is monotone, the proximal operator $P_k$ is single valued and well defined for all $x \in \dom{T}$. For an arbitrary maximal monotone operator the proximal point algorithm converges to an element of the set $T^{-1}(0)$ if it is nonempty. The proximal point algorithm also allows for approximate evaluation of $P_k$, a key consideration when designing a practical algorithm. It was shown in \cite{rockafellar1976monotone} that the proximal point algorithm can tolerate errors which satisfy the following
\begin{equation} \label{eq:prox_error}
  ||x_{k+1} - P_k(x_k)|| \leq \delta_k \sigma_k, \quad \sum_{k= 0}^\infty \delta_k <\infty.
\end{equation}
In general, if $T^{-1}(0) = \emptyset$ the algorithm diverges. However, we show in Section~\ref{ss:infeasibility} that for the special case considered here the proximal point algorithm can be used to detect infeasibility.

\section{The Stabilized Semismooth Algorithm}
In this section we describe the FBstab algorithm. The main idea is to regularize the original problem, solve it using a semismooth Newton-type method, then use the proximal point algorithm to iteratively refine the solution. The regularization ensures that each proximal subproblem has a unique primal-dual solution and satisfies the regularity conditions needed to ensure fast convergence of the inner Newton-type solver. In addition, semismooth Newton-type methods can be warmstarted and terminated early. As a result, each proximal subproblem can be solved approximately and warmstarted with the solution of the previous one. This makes FBstab very efficient, often requiring only one to two Newton iterations to solve each proximal subproblem.

The FBstab algorithm is summarized in Algorithms~\ref{algo:fbstab}- \ref{algo:check_infeas}. Algorithm~\ref{algo:fbstab} implements the proximal point algorithm and is discussed in Section~\ref{ss:outer_prox}. Algorithm~\ref{algo:eval_prox} evaluates the proximal operator, $P_k$, using a semismooth Newton's method as discussed in Section~\ref{ss:inner_newton} and Algorithm~\ref{algo:check_infeas} checks for infeasibility as discussed in Section~\ref{ss:infeasibility}.

\begin{algorithm}[h]
\caption{The FBstab algorithm}
\label{algo:fbstab}
\begin{algorithmic}[1]
\renewcommand{\algorithmicrequire}{\textbf{Inputs:}}
\renewcommand{\algorithmicensure}{\textbf{Outputs:}}
\Require $\sigma,\tau_r,\tau_a,\tau_d > 0$, $\kappa \in (0,1)$, Initial guess $x_0 = (z_0,\lambda_0,v_0)$
\Ensure Primal-dual solution, $x^*$, or infeasibility status and certificate $\Delta x^*$
\Procedure{\texttt{FBstab}}{}
\State $x \gets x_0$, $k \gets 0$
\State $\Delta x \gets \infty$
\State $\epsilon_0 \gets ||\pi(x_0)||$
\State $\delta \gets \min(\epsilon_0/\sigma,1)$
\Repeat
  \State $\delta \gets \min(\kappa \delta ,\epsilon/\delta)$ \label{line:delta_tight}
  \State $x^+\gets$ \Call{EvalProx}{$x,\delta \sigma,\sigma$}
  \State $\Delta x \gets x^+ - x$
  \State \Call{CheckFeasibility}{$\Delta x$}
  \State $x \gets x^+$
  \State $\epsilon \gets ||\pi(x)||$, $k\gets k+1$
\Until{$\epsilon \leq \epsilon_0 \tau_r + \tau_{a}$ or $||\Delta x|| \leq \tau_d$}
\State \textbf{Stop:} $x^* \gets x$ is optimal
\EndProcedure
\end{algorithmic}
\end{algorithm}

\begin{algorithm}[h]
\caption{Evaluate the proximal operator}
\label{algo:eval_prox}
\begin{algorithmic}[1]
\renewcommand{\algorithmicrequire}{\textbf{Inputs:}}
\renewcommand{\algorithmicensure}{\textbf{Outputs:}}
\Require $\beta\in (0,1), \eta \in (0,0.5)$
\Procedure{\texttt{EvalProx}}{$\bar{x},\eps,\sigma$}
\State $x\gets \bar{x}$
\Repeat
  \State Compute $V \in \bar{\partial}R(x,\bar{x},\sigma)$, see Algorithm~\ref{algo:compute_cdiff}
  \State Solve $V\Delta x = -R(x,\bar{x},\sigma)$ for $\Delta x$, see \eqref{eq:pfb_mapping}
  \State $t\gets 1$
  \While{$\theta(x+t\Delta x) \geq \theta(x) + \eta t \nabla \theta(x)$} \label{line:linesearch}, see \eqref{eq:merit}
    \State $t\gets \beta t$
  \EndWhile
  \State $x \gets x+ t \Delta x$
\Until{$||F(x,\bar{x},\sigma)||\leq \eps \min\{1,x-\bar{x}\}$} \label{line:dx_tight}
\State \Return x
\EndProcedure
\end{algorithmic}
\end{algorithm}

\begin{algorithm}[h]
\caption{Check for infeasibility}
\label{algo:check_infeas}
\begin{algorithmic}[1]
\renewcommand{\algorithmicrequire}{\textbf{Inputs:}}
\renewcommand{\algorithmicensure}{\textbf{Outputs:}}
\Require $\tau_{inf} > 0$
\Procedure{\texttt{CheckFeasibility}}{$\Delta x$}
\State $(\Delta z,\Delta v,\Delta v) \gets \Delta x$
\If{$||H \Delta z||_{\infty} \leq \tau ||\Delta z||_\infty$, $f^T\Delta z < 0$, \\ $\quad \quad \mathrm{max}(A\Delta z) \leq 0$, $||G\Delta z||_\infty \leq \tau ||\Delta z||_\infty$}
\State $\Delta x^* \gets (\Delta z,\Delta \lambda, \Delta v)$
\State \textbf{Stop:} $\Delta z^*$ certifies dual infeasibility
\EndIf
\If{$||A^T v + G^T \lambda||_\infty \leq \tau (||\Delta v|| + ||\Delta \lambda||)$, \\ $\quad \quad\Delta v^T b + \lambda^T h < 0$}
\State $\Delta x^* \gets (\Delta z,\Delta \lambda, \Delta v)$
\State \textbf{Stop:}  $(\Delta \lambda^*,\Delta v^*)$ certifies primal infeasibility
\EndIf
\EndProcedure
\end{algorithmic}
\end{algorithm}

\begin{rmk}
Both semismooth Newton methods \cite{fischer1992special} and proximal point algorithm \cite{rockafellar1976monotone} have been extensively studied in the literature. Our main contribution is the novel synergistic combination of the two methods with infeasibility detection techniques originally derived for ADMM \cite{banjac2017infeasibility} to produce a method that is very general, theoretically justified, and effective in practice. 
\end{rmk}

\begin{rmk}
This method could theoretically be extended to more general smooth convex programs. However, then the curvature of the constraint Hessians would enter into the subproblems and we would need to maintain positivity of the dual variables during the Newton iterations to ensure non-singularity of the generalized Jacobians. This would make warmstarting the algorithm difficult; as a result we have elected to focus on QPs. Moreover, due to the polyhedrality of the solution set, we are able to establish a stronger convergence rate results for QPs than for more general convex programs (see Theorem~\ref{thrm:FBstab_convergence}).
\end{rmk}

\subsection{Outer Proximal Point Iterations} \label{ss:outer_prox}

The KKT conditions of \eqref{eq:QP} can be rewritten as the following variational inequality (VI)
\begin{subequations} \label{eq:KKTVI}
\begin{gather}
\nabla_z L(z,\lambda,v) = 0,\\
h - Gz =0,\\
b - Az + N_+(v) \ni 0 
\end{gather}
\end{subequations}
where $N_+$ is the normal cone of the nonnegative orthant;  \eqref{eq:KKTVI} can be compactly expressed as,
\begin{equation} \label{eq:VI}
  T(x) = F(x) + N(x) \ni 0, 
\end{equation}
where $N$ is the normal cone of $\Gamma= \reals^n \times \reals^m \times \reals^q_{\geq0}$ and 
\begin{equation} \label{eq:base_def}
  F(x) = \begin{bmatrix}
    H & G^T & A^T\\-G & 0 & 0\\ -A & 0 & 0
  \end{bmatrix}\begin{bmatrix}
    z\\\lambda\\v
  \end{bmatrix} + \begin{bmatrix}
    f\\h\\b
  \end{bmatrix} = Kx + w .
\end{equation}
Solving \eqref{eq:QP} is equivalent to finding an element of $T^{-1}(0)$. 

\medskip
\begin{prp} \label{prp:T_properties}
The variational inequality \eqref{eq:KKTVI} has the following properties: (i) It is maximal monotone. (ii) If nonempty, its solution set, $T^{-1}(0)$, is closed and convex.
\end{prp}
\medskip
\begin{pf*}{Proof.}
(i): The base mapping \eqref{eq:base_def} is maximal monotone since it is single valued, affine and $K^T +K \succeq 0$, where $K$ is defined in \eqref{eq:base_def}. Variational inequalities of the form \eqref{eq:VI} are maximal monotone if the single valued portion $F$ is monotone  \cite{rockafellar1970maximality}. (ii): See \cite{minty1962monotone}.
\medskip
\end{pf*}

Since \eqref{eq:KKTVI} is maximally monotone we can apply the proximal point algorithm to \eqref{eq:VI}, this is the proximal methods of multipliers \cite{rockafellar1976augmented}. The proximal operator $P_k(x)$ can be evaluated by finding $x$ satisfying 
\begin{equation} \label{eq:prox_subproblem}
  T_\sigma(x) = F(x) + \sigma_k(x - x_k) + N(x) \ni 0,
\end{equation}
which is itself a variational inequality. However, due to the regularization term, \eqref{eq:prox_subproblem} is guaranteed to have a unique solution and to satisfy certain useful regularity properties (Proposition~\ref{prp:R_properties}). As a result, we can construct a semismooth Newton solver for the subproblems with a quadratic rate of convergence (Theorem~\ref{thm:newton_convergence}). In addition, the proximal point algorithm allows for approximate evaluation of $P_k$ and we warmstart the semismooth Newton solver at each iteration. Taken together, these measures allow FBstab to evaluate the proximal operator efficiently. The norm of the natural residual function, 
\begin{equation} \label{eq:nat_res}
  \pi(x) = x - \Pi_C(x - F(x)),
\end{equation}
 where $\Pi_C$ denotes euclidean projection onto the closed, convex set $C$, is used as a stopping criterion for the algorithm. It is a local error bound \cite[Theorem 18]{pang1997error}, i.e\footnote{For a closed set $C$, $\dist{x}{C} = \inf_{\bar{x}}\{x- \bar{x}~|~\bar{x}\in C\}$.}., $\dist{x}{T^{-1}(0)} = \Oh(||\pi(x)||)$.

\subsection{Inner Semismooth Newton Solver} \label{ss:inner_newton}

The core of FBstab is the inner solver which efficiently evaluates the proximal operator by applying a Newton-type method to a semismooth reformulation of \eqref{eq:prox_subproblem}. We construct an appropriate reformulation using a so-called nonlinear complementarity problem (NCP) function \cite{sun1999ncp}. A NCP function $\phi:\reals^2 \to \reals$ has the property that
\begin{equation}
   \phi(a,b) = 0 \quad \Leftrightarrow \quad a \geq 0,~b\geq 0,~ ab = 0.
\end{equation} 
In this paper we use the penalized Fischer-Burmeister (PFB) function \cite{chen2000penalized},
\begin{equation}
  \phi(a,b) = \alpha(a + b - \sqrt{a^2 + b^2}) + (1-\alpha)a_+ b_+,
\end{equation}
where $\alpha \in (0,1)$ is fixed and $x_+$ denotes projection onto the nonnegative orthant. The penalized FB (PFB) function is similar to the Fischer-Burmeister\cite{fischer1992special} function but has better theoretical and numerical properties \cite{chen2000penalized}. Using this NCP function we can construct the following mapping,
\begin{subequations} \label{eq:pfb_mapping}
\begin{gather} 
  R_k(x) = R(x,x_k,\sigma_k) = \begin{bmatrix}
    \nabla_zL(x) + \sigma_k(z-z_k)\\ h - Gz + \sigma_k (\lambda - \lambda_k)\\ \phi(y,v)
  \end{bmatrix},\\
  y = b-Az + \sigma_k(v - v_k)
\end{gather}
\end{subequations}
where the NCP function is applied elementwise. We will also make use of the following merit function
\begin{equation} \label{eq:merit}
  \theta_k(x) = \frac12||R_k(x)||_2^2.
\end{equation}
The properties of $R_k$ and $\theta_k$ are summarized below.
\medskip
\begin{prp} \label{prp:R_properties}
The function, $R_k$ in \eqref{eq:pfb_mapping}, and its merit function, $\theta_k(x)$, have the following properties:
\begin{enumerate}[1)]
  \item $R_k$ is strongly semismooth on $\reals^l$.
  \item $R(x_k^*,x_k,\sigma_k) = 0$ if and only if $x_k^* = P_k(x_k)$. Further, $x_k^*$ is unique and exists irrespective of the problem data.
  \item $||R_k(x)||$ is a global error bound, i.e., there exists $\tau > 0$ such that $||x - x_k^*|| \leq \tau||R_k(x)||$
  \item $\theta_k$ is continuously differentiable and, for any $V \in \bar{\partial} R_k(x)$, its gradient is $ \nabla\theta_k(x) = V^TR_k(x)$.
\end{enumerate}
\end{prp}
\medskip
\medskip
\begin{pf*}{Proof.}
\textit{1)} The PFB function is strongly semismooth \cite[Proposition 2.1]{chen2000penalized} and is composed with affine functions to form $R$. Strong semismoothness of $R$ then follows from the composition rules for semismooth functions, see, e.g., \cite[Propositions 1.73 and 1.74]{izmailov2014newton}. 

\textit{2)} The VI \eqref{eq:prox_subproblem} is defined by the sum of a monotone and strongly monotone operator and is thus strongly monotone. Strongly monotone operators always have a unique zero \cite{minty1962monotone}. The zeros of $R$ exactly coincide with those of \eqref{eq:prox_subproblem} by the properties of NCP functions.

\textit{3)} Let $\mathcal{K} = \{x~|~(K + \sigma I) x + N_\infty(x) \ni 0\}$ denote the ``kernel'' of \eqref{eq:prox_subproblem}, where $K$ is defined in \eqref{eq:base_def}, and $N_\infty(x)$ is the normal cone of $\Gamma_\infty$, the recession cone of $\Gamma$. Since $\Gamma$ is a convex cone, $\Gamma_\infty = \Gamma$. Theorem 20 of \cite{pang1997error} states that the norm of the natural residual function $\pi_k(x)$, defined in \eqref{eq:natres}, is a global error bound if $\mathcal{K} = \{0\}$. Its clear that $x = 0$ satisfies $(K + \sigma I) x + N_\infty(x) \ni 0$ and since $K + \sigma I$ is strongly monotone, the solution must be unique. Thus, applying \cite[Theorem 20]{pang1997error} there exists $\tau_1 > 0$ such that $||x-x_k^*|| \leq \tau_1 ||\pi_k(x)||$. The equivalence of $||\pi_k(x)||$ and $||R_k(x)||$, when used as an error bound, can then be established using the same arguments as \cite[Theorem 3.11]{chen2000penalized}; we omit the details for brevity.

\textit{4)} See \cite[Theorem 3.2]{chen2000penalized} or \cite[Proposition 2]{liaomcphersonFBRS}.
\medskip
\end{pf*}
\medskip
The inner solver evaluates $P_k$ by solving the rootfinding problem $R_k(x) = 0$ using a damped semismooth Newton's method. The inner iterative scheme\footnote{$i$ and $k$ are used for inner and outer iterations respectively.}, is
\begin{equation} \label{eq:ss_newton}
  x_{i+1|k} = x_{i|k} - t V^{-1}R_k(x_{i|k}), \quad V \in \bar{\partial} R_k(x_{i|k}),
\end{equation}
where $t \in (0,1]$ is a step length, chosen using a backtracking linesearch, which enforces global convergence. We will show in the next section that the matrix $V$ is always invertible and thus the iteration \eqref{ss:inner_newton} is well defined. In particular, all elements of $\cdiff R_k(x^*)$, where $x^*$ is the root, are non-singular which leads to quadratic convergence of \eqref{eq:ss_newton}. The rootfinding problem, $R_k(x) = 0$, has a unique solution even if the original QP is degenerate or infeasible due to the outer proximal point layer.
\begin{rmk} \label{rmk:natres}
The natural residual mapping
\begin{equation} \label{eq:natres}
   \pi_k(x) = x - \Pi_\Gamma(x - F(x) + \sigma_k(x-x_k)),
\end{equation} 
is similar to \eqref{eq:pfb_mapping} with $\phi(y,v)$ replaced by $\mathrm{min}(y,v)$ where the min operation is applied component wise.
\end{rmk}

\subsection{The Newton step system} 
In this section, we analyze the properties of the C-subdifferential and the associated Newton step system. We begin with the following proposition which establishes some properties of the C-subdifferential.
\begin{prp}
For any $x,\bar{x}\in \reals^l$, $\sigma >0$, any $V \in \cdiff R(x,\bar{x},\sigma)$ is of the form
\begin{equation} \label{eq:cdiff_estimate}
  V = \begin{bmatrix}
  H_\sigma & G^T & A^T \\
  -G & \sigma I & 0\\
  -CA& 0 & D
  \end{bmatrix},
\end{equation}
where $H_\sigma = H + \sigma I$, and $C = diag(\gamma_j)$, $D =diag(\mu_j + \sigma \gamma_j)$, are diagonal matrices with entries $(\gamma_j,\mu_j) \in \partial \phi(y_j,v_j)$. 
\end{prp}
\medskip
\begin{pf*}{Proof.}
The proof follows \cite[Proposition 2.3]{chen2000penalized}. By definition 
\begin{equation}
  \cdiff R(x,\bar{x},\sigma) = \partial R_{1}(x,\bar{x},\sigma) \times ~...~\times \partial R_l(x,\bar{x},\sigma),
\end{equation}
thus we need only to characterize the generalized gradients. The first two blocks are continuously differentiable so $\partial R_{i}(x,\bar{x},\sigma) = \{\nabla R_{i}(x,\bar{x},\sigma)^T\}$, $i = 1,\ldots, n+m$. The last block satisfies
\begin{equation}
  V_i= \begin{bmatrix}
    -\gamma_i(x) A_i & 0 & (\mu_i(x) + \sigma \gamma_i(x)) e_i
  \end{bmatrix},
\end{equation}
where $e_i$ are rows of identity, by \cite[Proposition 2.1 and Theorem 2.3.9]{clarke1990optimization}.
\end{pf*}
\medskip

Explicit expressions for the generalized gradient of $\phi(a,b)$ are given by \cite[Proposition 2.1]{chen2000penalized}
\begin{multline*}
\partial\phi(a,b) = (\gamma_i,\mu_i)= \\
  \begin{cases}
  \alpha (1- \frac{a}{r},1-\frac{b}{r}) + (1-\alpha)(b_+\partial a_+, a_+ \partial b_+) & \text{if}~~r \neq 0\\
  \alpha(1-\eta,1-\zeta) & \text{if} ~~ r = 0
  \end{cases}
\end{multline*}
where $r = \sqrt{a^2 + b^2}$, $\eta$ and $\zeta$ are arbitrary numbers satisfying $\eta^2 + \zeta^2 = 1$, and 
\begin{equation}
  \partial u_+ \in 
  \begin{cases}
  \{1\} & \text{if } u > 0\\
  [0,1] & \text{if } u = 0\\
  \{0\} & \text{if } u < 0.
  \end{cases}
\end{equation}
A procedure for computing an element of $\partial R(x,\bar{x},\sigma)$ is given in Algorithm~\ref{algo:compute_cdiff}.

\begin{algorithm}[H]
\caption{Compute Generalized Jacobian}
\label{algo:compute_cdiff}
\begin{algorithmic}[1]
\renewcommand{\algorithmicrequire}{\textbf{Inputs:}}
\renewcommand{\algorithmicensure}{\textbf{Outputs:}}
\Require Points $x,\bar{x}$, Tolerance $\zeta > 0$, $\sigma > 0$
\Ensure $V \in \reals^{l\times l}$
\Procedure{getJacobian}{$x, \bar{x},\sigma$}
\State $y \gets b - Ax + \sigma (x-\bar{x}) $
\State $r_i \gets \sqrt{y_i^2 + v_i^2}$
\State Define $S_0 = \{i~|~r_i \leq \zeta\}$
\State Define $S_1 = \{i~|~ y_i > 0, v_i > 0\}$
\For{$ i = 1,2, ... ~q$}
\If{$i \in S_0$}
\State $\gamma_i \gets \alpha (1- 1/\sqrt{2})$,~~$\mu_i \gets \alpha (1- 1/\sqrt{2})$
\ElsIf{$i \in S_1$}
\State $\gamma_i \gets \alpha (1- y_i/r_i) + (1-\alpha) v_i$,
\State $\mu_i \gets \alpha (1- v_i/r_i) + (1-\alpha) y_i$
\Else
\State $\gamma_i \gets \alpha (1- y_i/r_i)$,~~ $\mu_i \gets \alpha (1- v_i/r_i)$
\EndIf
\EndFor
\Statex
\State $V \gets \begin{bmatrix}
  H + \sigma I & G^T & A^T \\
  -G & \sigma I & 0\\
  -\diag{(\gamma)}A& 0 & \diag{(\mu + \sigma \gamma)}
  \end{bmatrix}$
\EndProcedure
\end{algorithmic}
\end{algorithm}

Next we consider the regularity properties of $V$ which are critical to the behaviour of the semismooth Newton solver. Due to the stabilizing effect of the outer proximal point algorithm all elements of $\bar{\partial}R(x,\bar{x},\sigma)$ are nonsingular, a property we will refer to as C-regularity.
\begin{dfn}
A mapping $G:\reals^N \to \reals^N$ is C-regular at a point $x\in \reals^N$ if it is semismooth at $x$ and all $V \in \cdiff G(x)$ are non-singular.
\end{dfn}

\begin{thm}[Regularity of the C-subdifferential] \label{thrm:Creg}
$R_k(x) = R(x,x_k,\sigma)$ is C-regular for any $x,x_k\in \reals^l$ and  $\sigma > 0$.
\end{thm}

\medskip
\begin{pf*}{Proof.}
For any $V\in \bar{\partial}R_k(x)$, the Newton step system \eqref{eq:ss_newton} has the form
\begin{equation} \label{eq:newton_expanded}
  \begin{bmatrix}
  H_\sigma & G^T & A^T \\
  -G & \sigma I & 0\\
  -CA& 0 & D
  \end{bmatrix}
  \begin{bmatrix}
    \Delta z\\ \Delta \lambda \\ \Delta v
  \end{bmatrix} = 
  \begin{bmatrix}
    r_1\\r_2\\r_3
  \end{bmatrix},
\end{equation}
where $R_k(x) = -[r_1^T~~r_2^T~~r_3^T]^T$ as in \eqref{eq:pfb_mapping}. For all $j \in \{1, ... q\}$ we have that $\mu_j \geq 0,~\gamma_j \geq 0$ and $(\mu_j,\gamma_j) \neq 0$. Thus $D_{jj} = \mu_j + \sigma \gamma_j > 0$ implying $D \succ 0$. Since $D \succ 0$ we can eliminate the third row of \eqref{eq:newton_expanded} algebraically and negate the second leading to the following pair of linear systems of equations
\begin{subequations} \label{eq:newton_reduced}
\begin{gather} 
  \begin{bmatrix}
  E & G^T\\
  G & -\sigma I \\
  \end{bmatrix}
  \begin{bmatrix}
    \Delta z\\\Delta \lambda
  \end{bmatrix} = \begin{bmatrix}
    r_1 - A^TD^{-1}r_3 \\ -r_2
  \end{bmatrix},\\
  D\Delta v = r_3 + CA\Delta z.
\end{gather}
\end{subequations}
The matrices $E = H_\sigma + A^T CD^{-1}A$ and $\sigma I$ are positive definite so the block $2 \times 2$ matrix in \eqref{eq:newton_reduced} is symmetric quasidefinite and thus invertible \cite{vanderbei1995symmetric}. As a result, \eqref{eq:newton_expanded} has a unique solution, implying that $V$ is nonsingular.
\medskip
\end{pf*}
The C-regularity property of $R_k$ guarantees that the inner iterations are globally well defined and that the iteration \eqref{eq:ss_newton} will converge at a quadratic rate to the unique solution of \eqref{eq:prox_subproblem} (Theorem~\ref{thm:newton_convergence}). 

\section{Convergence Analysis} \label{ss:convergence}
In this section we discuss the convergence of the FBstab algorithm. First, we prove convergence of the inner Newton-type solver. Our analysis follows the established convergence theory for damped semismooth Newton's methods, see e.g., \cite{qi1996c,chen2000penalized,facchinei1997new} and is included for completeness.
\begin{thm}[Inner solver convergence]  \label{thm:newton_convergence}
Consider an arbitrary but fixed iteration $k$ of the outer proximal point algorithm. Suppose $x_0 = x_{0|k} = x_k \in \reals^l$ and let the sequence $\{x_i\}$ be generated by the \texttt{EvalProx} procedure in Algorithm~\ref{algo:fbstab}. Then:
\begin{enumerate}[i.]
  \item The sequence $\{x_i\}$ is well defined and converges to the unique point $x^*$ satisfying $x^* = P_k(x_0)$.
  \item The asymptotic rate of convergence is quadratic i.e.,
  \begin{equation*}
   	||x^* - x_{i+1}|| = \Oh(||x^* - x_i||^2) \text{ as } i \to \infty.
  \end{equation*} 
\end{enumerate}
\end{thm}
\medskip
\begin{pf*}{Proof.} Suppose that \texttt{EvalProx} generates an infinite sequence $\{x_i\}$ satisfying $||R_k(x_i)|| > 0$; if $R_k(x_i) = 0$ the algorithm will terminate. Let $x^*$ denote the unique solution of $\eqref{eq:prox_subproblem}$. The sequence $\{x_i\}$ is uniquely defined because $R_k$ is globally C-regular by Theorem~\ref{thrm:Creg}.

We begin with a local analysis. Suppose for now that $t = 1$, define $e_i = x_i - x^*$, and consider the update equation
\begin{subequations}
\begin{align*}
  ||e_{i+1}|| &= ||x_i - x^* - V^{-1}R_k(x_i)||, \\
  \leq &||V^{-1}||~||V(x_i - x^*) - R_k(x_i) + R_k(x^*)||,\\
  \leq& ||V^{-1}|| \sqrt{\sum_{j=1}^l ||V_j (x_i-x^*) - R_{k,j}(x) + R_{k,j}(x^*)||^2},
\end{align*}
where we have used that $\Delta x_i = -V^{-1}R_k(x_i)$ and the boundedness of $||V^{-1}||$, which holds by the global C-regularity of $R_k$, see Theorem~\ref{thrm:Creg}. Due to the strong semismoothness of $R_k$ and the properties of the C-subdifferential \cite{qi1996c} we have that,
\begin{equation*}
  ||V_j (x_i-x^*) - R_{k,j}(x) + R_{k,j}(x^*)|| = \Oh(||x_i - x^*||^2),
\end{equation*}
as $x_i \to x^*$. Thus we have that,
\begin{equation} \label{eq:Q_convergence}
  ||x_{i+1} - x^*|| = \Oh(||x_i - x^*||^2),
\end{equation}
as $x_i \to x^*$; this establishes local quadratic convergence. 

Next we consider global convergence. Since $\theta_k$ is continuously differentiable and $\nabla \theta_k(x_i) = V^T R_k(x_i)$, see Proposition~\ref{prp:R_properties}, we have that
\begin{align} \label{eq:descent_direction}
  \nabla \theta_k(x_i)^T \Delta x_i &= - (V^T R_k(x_i))^TV^{-1} R_k(x_i) \\
   & = -||R_k(x_i)||^2 \leq 0.
\end{align}
Performing a Taylor expansion yields,
\begin{equation}
  \theta_k(x_i + t\Delta x_i)= \theta_k(x_i) + t\nabla\theta(x_i)^T\Delta x_i + \oh(t).
\end{equation}
Since $\nabla\theta(x_i)^T\Delta x_i < 0$ by \eqref{eq:descent_direction} this implies that the linesearch condition (Line ~\ref{line:linesearch}, Algorithm~\ref{algo:eval_prox})
\begin{equation} \label{eq:linsesearch}
  \theta_k(x_i + t\Delta x_i) \geq \theta(x_i) + t \eta \nabla\theta(x_i)^T\Delta x_i
\end{equation}
must hold as $t \to 0$ ensuring that the linesearch terminates with $t > 0$. Note that
\begin{multline} \label{eq:armijo}
  \theta_k(x_i) + t \eta \nabla\theta_k(x_i)^T\Delta x_i =  \\
    \frac12 ||R_k(x_i)||^2 - \eta t ||R_k(x_i)||^2 
   = (1-2\eta t) \theta_k(x_i) 
\end{multline}
where we have used \eqref{eq:descent_direction}. Combining \eqref{eq:linsesearch} with \eqref{eq:armijo} and noting that $t \in (0,1]$, we have that
\begin{align}
  \theta_k(x_i + t\Delta x_i) &\leq  (1-2\eta t) \theta_k(x_i) < \theta_k(x_i),
\end{align}
since $\eta \in (0,0.5)$. As the sequence $\{\theta_k(x_i)\}$ is strictly decreasing and bounded below by zero it converges to some $\theta^* \geq 0$ and $\lim_{i\to \infty} \theta_k(x_{i+1})- \theta_k(x_i) = 0$. Since
\begin{equation} \label{eq:global1}
  \theta_k(x_{i+1}) - \theta_k(x_i) \leq t \eta \nabla \theta_k(x_i)^T \Delta x_i,
\end{equation}
we have, rearranging \eqref{eq:global1} and using \eqref{eq:descent_direction}, that
\begin{subequations}
\begin{align} \label{eq:global_conv}
  \theta_k(x_i) &= \frac12 ||R_k(x_i)||^2,\\
  & = -\frac12 \nabla \theta_k(x_i)^T \Delta x_i,\\
  &\leq (\theta_k(x_i) - \theta_k(x_{i+1}))/t_i \eta \to 0,
\end{align}
\end{subequations}
as $i \to \infty$. Since $\theta_k(x) = 0$ if and only if $x = x^*$ and $\theta_k$ is continuous this establishes global convergence.

Finally, it remains to show that the algorithm will recover the local quadratic converge rate established earlier. Since $\{x_i\} \to x^*$ as $i \to \infty$ eventually $x_i$ will be sufficiently close to $x^*$ for \eqref{eq:Q_convergence} to hold if a unit step is admitted by the linesearch. Algorithm~\ref{algo:eval_prox} always begins by checking $t = 1$. Since,
\begin{align*}
\theta_k(x_i+\Delta x_i) &= \frac12 ||R_k(x_i + \Delta x_i) - R_k(x^*)||^2,\\
& = \Oh(||x_i + \Delta x_i - x^*||^2),\\
& = \Oh(||x_i -x^*||^4) = \oh(||x_i -x^*||^2),\\
& = \oh(||R_k(x_i)||^2) = \oh(\theta_k(x_i)),
\end{align*}
where the second line holds by Lipschitz continuity of $R_k$,  the third by \eqref{eq:Q_convergence} and the fourth by the error bound $||x_i - x^*|| = \oh(||R_k(x_i)||)$ in Proposition~\ref{prp:R_properties}, the linesearch will eventually accept unit steps implying that \eqref{eq:Q_convergence} holds as $i \to \infty$.
\end{subequations}
\medskip
\end{pf*}

The following theorem summarizes the convergence properties of FBstab when $T^{-1}(0) \neq \emptyset$, i.e., when \eqref{eq:QP} has a primal-dual solution. We analyze the behaviour of FBstab when $T^{-1}(0) = \emptyset$ in Section~\ref{ss:infeasibility}.
\medskip
\begin{thm}[Convergence of FBstab] \label{thrm:FBstab_convergence} 
Let $x_0 \in \reals^l$ be arbitrary, suppose $T^{-1}(0) = (F + N)^{-1}(0)$ is nonempty and let $\{x_k\}$ be generated by FBstab. Then $\{x_k\} \to x^* \in T^{-1}(0)$ as $k\to \infty$. Further, the convergence rate is at least linear and if, in addition, $\sigma_k \to 0$ as $k\to \infty$ then the convergence rate is superlinear.
\end{thm}
\medskip
\begin{pf*}{Proof.} 
FBstab is an instance of the proximal point algorithm so we can employ\cite[Theorem 2.1]{luque1984asymptotic} to establish both convergence and the convergence rate. The error bound condition needed by \cite[Theorem 2.1]{luque1984asymptotic} ($(A_r')$ in \cite{luque1984asymptotic}) is
\begin{equation} \label{eq:AR}
	\dist{0}{T_\sigma(x_{k+1})} \leq \delta_k \sigma_k \min\{1,||x_{k+1}-x_k||\}
\end{equation}
where $T_\sigma$ is defined in \eqref{eq:prox_subproblem}, and $\sum_{k=0}^{\infty} \delta_k < \infty$. The inequality \eqref{eq:AR} is enforced by construction (Line~\ref{line:dx_tight} of Algorithm~\ref{algo:eval_prox}) since $||R_k(x)||$ is an error bound for the subproblems by Prop~\ref{prp:R_properties}. Moreover, ${\delta_k}$ satisfies, $\delta_k \leq {\kappa} ^k$ (Line~\ref{line:delta_tight} of Algorithm~\ref{algo:fbstab}) thus $\kappa \in (0,1)$, implies $\sum_0^{\infty} \delta_k = \frac{1}{1-\kappa} < \infty$. It remains to show that there exists $a,\eps> 0$ such that for all $ w \in \eps \mathbb{B}$ 
\begin{equation}
\dist{x}{T^{-1}(0)} \leq a ||w||, ~\forall x\in T^{-1}(w), 
\end{equation}
or equivalently (see e.g., \cite[Section 3D]{dontchev2009implicit})
\begin{equation}
	T^{-1}(w) \subseteq T^{-1}(0) + a ||w|| \mathbb{B},
\end{equation}
where $\mathbb{B} = \{x~|~||x|| \leq 1\}$. This property actually holds globally because $T(x) = F(x) + N(x)$ is a polyhedral variational inequality, see \cite[Section 3D]{dontchev2009implicit}. Thus we can invoke \cite[Theorem 2.1]{luque1984asymptotic} to conclude that FBstab converges, the rate of convergence is globally linear, and that the convergence rate is superlinear if $\sigma_k \to 0$ as $k \to \infty$.
\end{pf*}
\medskip

Finally, we state the following theorem which provides rigorous justification for the subproblem warmstarting strategy employed in FBstab.
\medskip
\begin{thm}[Lipschitz continuity of subproblems] \label{thm:subproblem_regularity}
Let the sequence $\{x_k\}$ be generated by FBstab with $x_0$ arbitrary. Then the proximal operator is Lipschitz continuous, i.e., at any iteration $k$ the proximal operator $P_k$ satisfies
\begin{gather}
  ||P_k(x_k) - P_{k-1}(x_{k-1})|| \leq \eta ||x_k - x_{k-1}||, \\
  \eta^{-1} = \lambda_{min}\left(\frac{K^T + K + 2\sigma I}{2\sigma}\right),
\end{gather}
where $\lambda_{min}(\cdot)$ designates the smallest eigenvalue of a symmetric matrix and $K$ is defined in \eqref{eq:base_def}.
\end{thm}
\medskip
\begin{pf*}{Proof.}
The variational inequality \eqref{eq:prox_subproblem} can be written as $F(x)/\sigma + x +N(x) \ni x_k$ which is a parameterized variational inequality with $x_k$ as the parameter. Its strong monotonicity constant is $\eta > 0$. The result then follows from \cite[Theorem 2F.6]{dontchev2009implicit}.
\end{pf*}
\medskip
Each proximal subproblem computes $P_k(x_k)$ by solving $R_k(x) = 0$ starting from $x_k$ as an initial guess. Theorem~\ref{thm:subproblem_regularity} implies that eventually $||P_k(x_k) - x_k||$ will become sufficiently small so that the quadratic convergence rate of Theorem~\ref{thm:newton_convergence} holds immediately and the semismooth method converges rapidly. This happens because,
\begin{align*}
  &||P_k(x_k) - x_k||,\\ &
  = ||[P_k(x_k) - P_{k-1}(x_{k-1})] - [x_k - P_{k-1}(x_{k-1})]||,\\
  & \leq ||x_k - P_{k-1}(x_{k-1})|| + ||P_{k}(x_{k}) - P_{k-1}(x_{k-1})||,\\
  & \leq ||x_k - P_{k-1}(x_{k-1})|| + \eta ||x_{k-1} - x_{k-2}||,
\end{align*} and, since the algorithm is converging, $||x_{k-1} - x_{k-2}|| \to 0$ and $||x_k - P_{k-1}(x_{k-1})|| \to 0$ as $k\to \infty$. We observe this behaviour in practice, typically after the first or second proximal iteration each subsequent proximal subproblem takes only one or two Newton iterations to converge.

\subsection{Infeasibility Detection} \label{ss:infeasibility}
In this section we apply the techniques developed in \cite{banjac2017infeasibility} to characterize the behaviour of FBstab in the case where either the QP \eqref{eq:QP} or its dual \eqref{eq:dual_QP} are infeasible. These results are not specific to FBstab and hold whenever the proximal point algorithm is used to solve \eqref{eq:KKTVI}. We begin by recalling infeasibility conditions for \eqref{eq:QP} and \eqref{eq:dual_QP}.

\medskip
\begin{prp} \label{prp:infeas_conditions} (Infeasibility conditions)\\
Dual infeasibility: Suppose there exists a vector $z \in \reals^n$ satisfying $Hz = 0$, $Az \leq 0$, $Gz = 0$, and $f^T z <0$. Then \eqref{eq:dual_QP} is infeasible.\\
Primal infeasibility: Suppose there exits a vector $(\lambda,v)$ such that $G^T \lambda + A^T v = 0$ and $\lambda^Th + v_+^T b < 0$, where $v_+$ is the projection of $v$ onto the nonnegative orthant. Then the feasible set of \eqref{eq:QP} is empty.
\end{prp}
\medskip
\begin{pf*}{Proof.}
See e.g., \cite[Proposition 1]{banjac2017infeasibility}.
\end{pf*}
\medskip
Any vector satisfying the conditions of Proposition~\ref{prp:infeas_conditions} is a certificate of primal or dual infeasibility. We will show that the proximal point algorithm generates these certificates when appropriate. When \eqref{eq:QP} is feasible, dual infeasibility is the same as the primal problem being unbounded below. First we review the limiting behaviour of the proximal point algorithm. The following lemma summarizes some results for averaged nonexpansive operators, a class which includes the proximal operator.

\medskip
\begin{lmm}\label{lemma:op_convergence}
Let $T: \mathcal{D} \mapsto \mathcal{D}$ be an averaged nonexpansive operator. In addition, suppose $x_k$ is generated by $x_k = T^k(x_0)$, $x_0 \in \mathcal{D}$, define $\delta x_k = x_{k+1} - x_k$, and let $\delta x $  be the projection of $0$ onto $\mathrm{cl}~\mathrm{range}~(T-I_d)$ where $\mathrm{cl}$ denotes the closure of a set and $I_d$ denotes the identity operator. Then as $k \to \infty$ we have that:
\begin{enumerate}[i.]
\item $\frac1k x_k \to \delta x$
\item $\delta x_k \to \delta x$
\item If $\mathrm{Fix}~T \neq \emptyset$ then $x_k \to x^* \in \mathrm{Fix}~T$, where $\mathrm{Fix}~T$ denotes the fixed points of $T$.
\end{enumerate}
\end{lmm}
\medskip
\begin{pf*}{Proof.}
(i): \cite[Corrolary 2]{pazy1971asymptotic}. (ii), (iii): \cite[Fact 3.2]{bauschke2004finding} .
\medskip
\end{pf*}
An immediate corollary of this is that $\delta x$ in Lemma~\ref{lemma:op_convergence} satisfies $\delta x = 0$ if $\mathrm{Fix}~T \neq \emptyset$. The following proposition applies Lemma~\ref{lemma:op_convergence} to our specific situation.

\begin{prp} \label{prp:infeas}
Let the sequence $\{x_k\} = \{(z_k,\lambda_k,v_k)\}$ be generated by the proximal point algorithm and define $\delta x_k = x_k - x_{k-1}$. Then there exists $\delta x = (\delta z,\delta \lambda,\delta v) \in \reals^l$ such that $(\delta z_k,\delta \lambda_k,\delta v_k) \to (\delta z,\delta \lambda,\delta v)$ as $k\to \infty$ and also satisfies the following properties:
\begin{enumerate}[(i)]
  \item $H\delta z = 0$,
  \item $A\delta z \leq 0$,
  \item $G\delta z = 0$,
  \item $\delta v \geq 0$,
  \item $f^T \delta z = -\sigma ||\delta z|| \leq 0$,
  \item $\delta \lambda^T h + \delta v^T b \leq 0$,
  \item $G^T \delta \lambda + A^T \delta v = 0$.
\end{enumerate}
\end{prp}
\medskip
\begin{pf*}{Proof.}
The proximal operator is firmly non-expansive \cite{rockafellar1976monotone} and thus averaged, see e.g., \cite[rmk 4.34]{bauschke2017convex}. The convergence of $\delta x_k$ to $\delta x$ as $k\to \infty$ then follows from Lemma~\ref{lemma:op_convergence}. Note that $\lim_{k\to\infty}~ \frac1k \delta x_k = 0$ and $\lim_{k\to\infty}~ \frac1k x_k = \delta x$ which we will use often in the sequel. We begin by rewriting \eqref{eq:prox_subproblem} in the following form:
\begin{subequations} \label{eq:infeas}
\begin{gather}
  Hz_k + f + G^T \lambda_k + A^T v_k + \sigma \delta z_k = 0, \label{eq:infeas1}\\
  h - Gz_k + \sigma \delta \lambda_k = 0, \label{eq:infeas2}\\
  \langle b - Az_k + \sigma \delta v_k,v_k \rangle = 0,\label{eq:infeas3}\\
  v_k \geq0,~ b - Az_k + \sigma \delta v_k \geq 0. \label{eq:infeas4}
\end{gather}
\end{subequations}
Further, \eqref{eq:infeas} is satisfied exactly in the limit since the condition $\eps_k \to 0$ as $k \to \infty$ is enforced by construction in Algorithm~\ref{algo:fbstab}. We now proceed point by point. \\
(i): Taking inner products, multiplying \eqref{eq:infeas3} and \eqref{eq:infeas2} by $1/k$, taking the limit, and applying Lemma~\ref{lemma:op_convergence} yields
\begin{subequations} \label{eq:i1}
\begin{gather}
  \lim_{k\to\infty}~ \frac1k \langle b - Az_k + \sigma \delta v_k,v_k \rangle = \langle -A\delta z,\delta v \rangle = 0,\\
   \Rightarrow \delta v^T A \delta z = 0,\\
  \lim_{k\to\infty}~ \frac1k \langle h - Gz_k + \sigma \delta \lambda_k, \lambda_k \rangle = \langle -G\delta z, \delta \lambda \rangle = 0.\\ 
  \Rightarrow \delta \lambda^T G \delta z = 0.
\end{gather}
\end{subequations}
The same procedure applied to \eqref{eq:infeas1} yields
\begin{multline}
  \lim_{k\to\infty}~ \frac1k \langle Hz_k + f + G^T \lambda_k + A^T v_k + \sigma \delta z_k, \delta z \rangle  \\
   = \delta z^T H \delta z + \delta \lambda^T G \delta z + \delta v^T A \delta z = 0,
\end{multline}
combining this with \eqref{eq:i1} we obtain that, since $H \succeq 0$,
\begin{equation}
  \delta z^T H \delta z = 0 \Rightarrow H \delta z = 0.
\end{equation}
(ii): Multiplying the second inequality in \eqref{eq:infeas4} by $1/k$ and taking the limit yields
\begin{equation}
  \lim_{k\to\infty}~ \frac1k (b - Az_k + \sigma \delta v_k) = -A\delta z \geq 0 \Rightarrow A\delta z \leq 0.
\end{equation}
(iii): Multiplying \eqref{eq:infeas2} by $1/k$ and taking the limit yields
\begin{equation}
  \lim_{k\to\infty}~ \frac1k (h - Gz_k + \sigma \delta \lambda_k) = -G\delta z = 0 \Rightarrow G\delta z = 0.
\end{equation}
(iv): Multiplying \eqref{eq:infeas4} by $1/k$ and taking the limit yields
\begin{equation}
  \lim_{k\to\infty}~ \frac1k v_k = \delta v \geq 0.
\end{equation}
(v): Taking the inner product of \eqref{eq:infeas1} with $\delta z_k$ then taking the limit and applying \eqref{eq:i1} and (i) yields:
\begin{multline}
  \lim_{k\to\infty}~ \langle Hz_k + f + G^T \lambda_k + A^T v_k + \sigma \delta z_k, \delta z_k \rangle \\
  = \delta z^T H \delta z + \delta \lambda^TG\delta z + \delta v^T A\delta z + f^T \delta z + \sigma ||\delta z||^2_2\\ 
  = f^T \delta z + \sigma ||\delta z||^2_2 = 0 ~ \Rightarrow f^T \delta z = -\sigma ||\delta z||^2_2 \leq 0.
\end{multline}
(vi): Taking the inner product of \eqref{eq:infeas3} and $\delta v_k$, and taking the limit yields
\begin{multline}
  \lim_{k\to\infty}~\langle b - Az_k + \sigma \delta z_k,\delta v_k \rangle\\
  = b^T \delta v - \delta v^T A \delta z + \sigma||\delta v||^2_2 = 0,
\end{multline}
since $\delta v^TA\delta z = 0$ we have that
\begin{equation} \label{eq:i2}
  b^T \delta v = -\sigma ||\delta v||_2^2 \leq 0.
\end{equation}
Applying the same procedure to \eqref{eq:infeas2} yields,
\begin{multline}
  \lim_{k\to\infty}~\langle h - Gz_k + \sigma \delta \lambda_k,\delta \lambda_k \rangle \\
  = h^T \delta \lambda - \delta \lambda G \delta z + \sigma ||\delta \lambda||^2_2 = 0,
\end{multline}
since $\delta \lambda^TG\delta z = 0$ this implies $h^T \delta \lambda = -\sigma||\delta \lambda||^2_2$, combining this with \eqref{eq:i2} yields
\begin{equation}
  h^T \delta \lambda + b^T \delta v = -\sigma(||\delta \lambda||_2^2 + ||\delta v||_2^2) \leq 0.
\end{equation}
(vii): Dividing \eqref{eq:infeas1} by $k$ and taking the limit yields
\begin{multline}
  \lim_{k\to\infty}~ \frac1k(Hz_k + f + G^T \lambda_k + A^T v_k + \sigma \delta z_k) \\
  = H\delta z + G^T\delta \lambda + A^T \delta v  = 0,
\end{multline}
applying (i) we have that $H\delta z = 0$ so we obtain
\begin{equation}
  G^T \delta \lambda + A^T \delta v = 0,
\end{equation}
which completes the proof.
\end{pf*}
\medskip

Armed with Proposition~\ref{prp:infeas} we can prove the following theorem summarizing the behaviour of FBstab when \eqref{eq:QP} or \eqref{eq:dual_QP} is infeasible.
\medskip
\begin{thm}[Infeasibility Detection] \label{thrm:infeas_detection} 
Suppose that \eqref{eq:QP} is primal-dual infeasible, i.e., the solution set of \eqref{eq:VI} is empty. Suppose $x_0 \in \reals^l$ is arbitrary, let the sequence of iterates $\{x_k\} = \{z_k,\lambda_k,v_k\}$ be generated by FBstab, and define $\delta x_k = x_{k+1} - x_k$. Then $\delta x_k \to \delta x$  as $k \to \infty$ where $\delta x = (\delta z,\delta \lambda,\delta v)$ satisfies the following properties:
\begin{enumerate}[(i)]
\item If $\delta z \neq 0$ then the dual QP \eqref{eq:dual_QP} is infeasible and $\delta z$ satisfies the dual infeasibility conditions in Proposition~\ref{prp:infeas_conditions}. 
\item If $(\delta \lambda, \delta v) \neq 0$ then the primal QP $\eqref{eq:QP}$ is infeasible and $(\delta \lambda, \delta v)$ satisfies the primal infeasibility conditions in Proposition~\ref{prp:infeas_conditions}.
\item If $\delta x \neq 0$ and $(\delta \lambda, \delta v) \neq 0$ then \eqref{eq:QP} and \eqref{eq:dual_QP} are infeasible.
\end{enumerate}
\end{thm}
\medskip
\begin{pf*}{Proof.}
(i): Follows from points (i), (ii), (iii), and (v), of Proposition~\ref{prp:infeas}. Note that if $\delta z \neq 0$ then $f^T \delta z = -\sigma ||\delta z|| < 0$.

(ii): Follows from points (vi), (vii), and (iv) of Proposition~\ref{prp:infeas}. Note that since $\delta v > 0$ due to point (iv) of Proposition~\ref{prp:infeas}, the condition $b^T \delta v_+ + h^T \delta \lambda < 0$ simplifies to $b^T \delta v + h^T \delta \lambda$.

(iii): Follows from points (i) and (ii) above.
\end{pf*}
\medskip
Theorem~\ref{thrm:infeas_detection} justifies Algorithm~\ref{algo:check_infeas}. This feature allows FBstab to exit gracefully if there is no primal-dual solution. Infeasibility detection is also important in many applications e.g., in branch and bound algorithms for mixed integer QPs \cite{fletcher1998numerical}.

\section{Numerical Experiments}
In this section we illustrate the performance of FBstab with some numerical experiments. We solve instances of the following optimal control problem (OCP),
\begin{subequations} \label{eq:MPCQP}
\begin{gather}
\underset{x,u}{\mathrm{min.}}~\sum_{i=0}^{N} \frac12 \begin{bmatrix}
	x_i \\ u_i
\end{bmatrix}^T \begin{bmatrix}
	Q_i & S_i^T\\
	S_i & R_i
\end{bmatrix} \begin{bmatrix}
	x_i \\u_i 
\end{bmatrix} +
\begin{bmatrix}
	q_i \\ r_i
\end{bmatrix}^T \begin{bmatrix}
	x_i \\u_i 
\end{bmatrix},\\
\mathrm{s.t.} ~~ x_0 = \xi,\\
x_{i+1} = A_i x_i + B_i u_i + c_i, ~ i \in \intrng{0}{N-1},\\
E_i x_i + L_i u_i + d_i \leq 0,~ i \in \intrng{0}{N},
\end{gather}
\end{subequations}
where $A_i,Q_i \in \reals^{n_x \times n_x}$, $B_i, S_i \in \reals^{n_x \times n_u}$, $R_i \in \reals^{n_u \times n_u}$, $q_i,c_i \in \reals^{n_x}$, $r_i \in \reals^{n_u}$, $E_i \in \reals^{n_c \times n_x}$, $L_i \in \reals^{n_c \times n_u}$, $d_i \in \reals^{n_c}$, $x_i,\xi \in \reals^{n_x}$, $u_i \in \reals^{n_u}$, $x = (x_0, ..., x_N)$, and $u = (u_0, \dots ,u_N)$. We require that
\begin{equation}
	\begin{bmatrix}
		Q_i & S_i^T\\
		S_i & R_i
	\end{bmatrix} \succeq 0~~ \forall i \in \intrng{0}{N},
\end{equation}
so the problem is convex. This QP is large but sparse and is often called the simultaneous or multiple shooting form of the MPC problem \cite{rawlings2009model}. The QP is also often solved in the so-called condensed form,
\begin{subequations}\label{eq:condensed_QP}
\begin{gather} 
\underset{u}{\mathrm{min.}} \quad \frac12 u^T H u + f(\xi)^T u, \\
\mathrm{s.t} \quad Au \leq b(\xi),
\end{gather}
\end{subequations}
which is in the control variables only and can be derived by eliminating the state variables $x$ in \eqref{eq:MPCQP} using the dynamic equations see e.g., \cite[Section 2.3]{kouvaritakis2016model} or \cite{borrelli2017predictive}. We consider three linear MPC benchmark problems; their properties are summarized in Table~\ref{tab:problem_sizes}.

\begin{figure}[htbp]
	\centering
	\includegraphics[width=0.95\columnwidth]{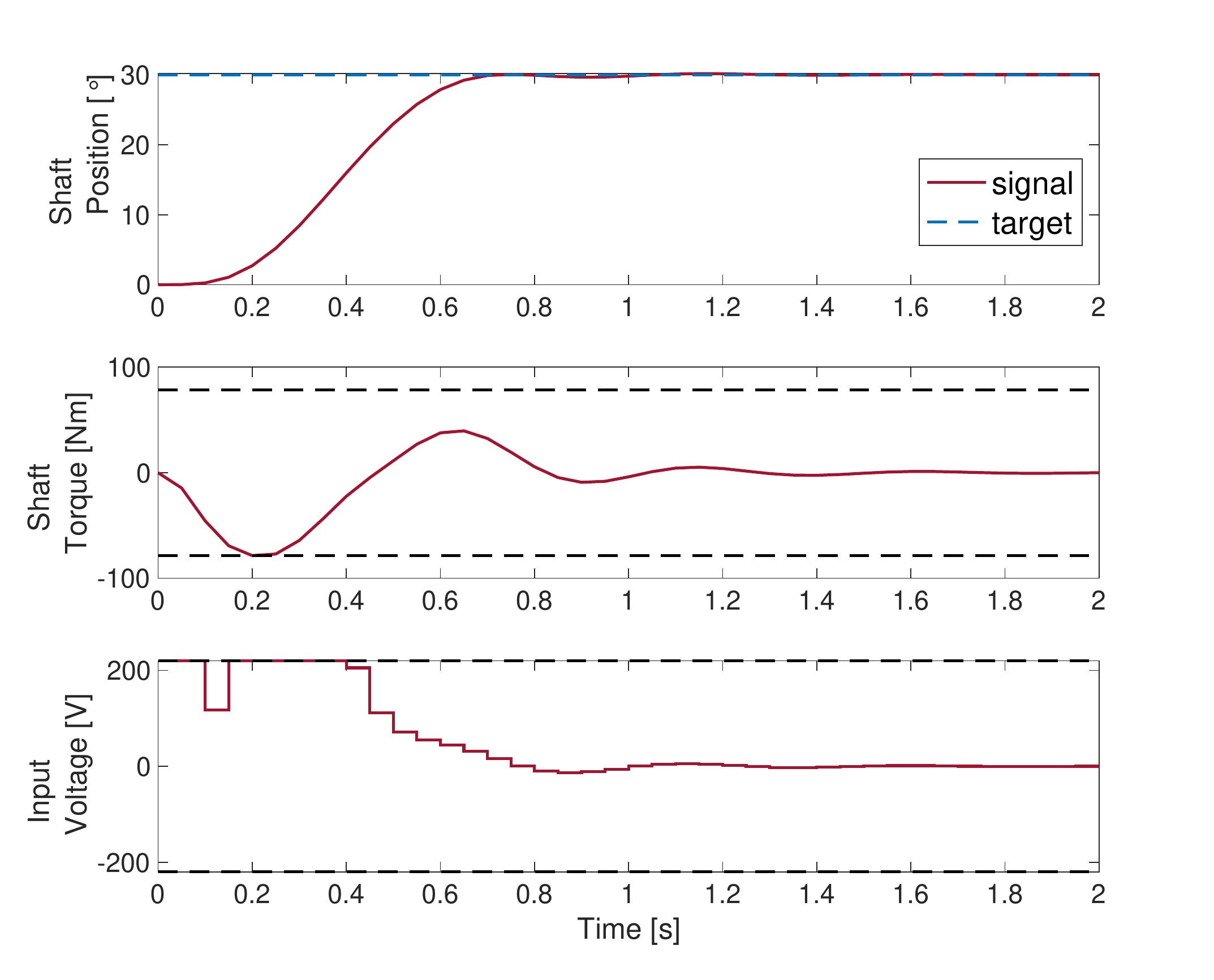}
	\caption{Closed-loop response in the servo motor example.}
	\label{fig:servo_cl}
\end{figure}

\textbf{Control of a Servo Motor\cite{bemporad1998fulfilling}:}
The objective is to drive the motor position $y_1$ to a desired angular position $r = 30^\circ$ while respecting the constraint $|y_{2,k}| \leq 78.5 ~Nm$ on the shaft torque and the constraint $|u| \leq 220~V$ on the motor input voltage. The continuous time model is
\begin{gather*}
\frac{d}{dt}x(t) = \begin{bmatrix}
0 & 1 & 0 & 1\\
-128 & -2.5 & 6.4 & 0\\
0 & 0 & 0 & 1\\
128 & 0 & -6.4 & -10.2
\end{bmatrix} x(t) + \begin{bmatrix} 
0\\0\\0\\1
\end{bmatrix}u(t) ,~~ \\ y(t) = \begin{bmatrix}
1 & 0 & 0 & 0\\
1282 & 0 & -64.0 & 0
\end{bmatrix}x(t) ,
\end{gather*}
which is discretized at $0.05$ s using a zero-order hold. The tuning matrices and initial condition are $Q_i = diag([10^3,0,0,0])$, $R_i = 10^{-4}$, and $x_0 = 0$. The traces of this model in closed-loop with an MPC controller are shown in Figure~\ref{fig:servo_cl}, the shaft angular position is driven to the reference while respecting the constraints on the shaft torque and input voltage.

\begin{figure}[htbp]
	\centering
	\includegraphics[width=0.95\columnwidth]{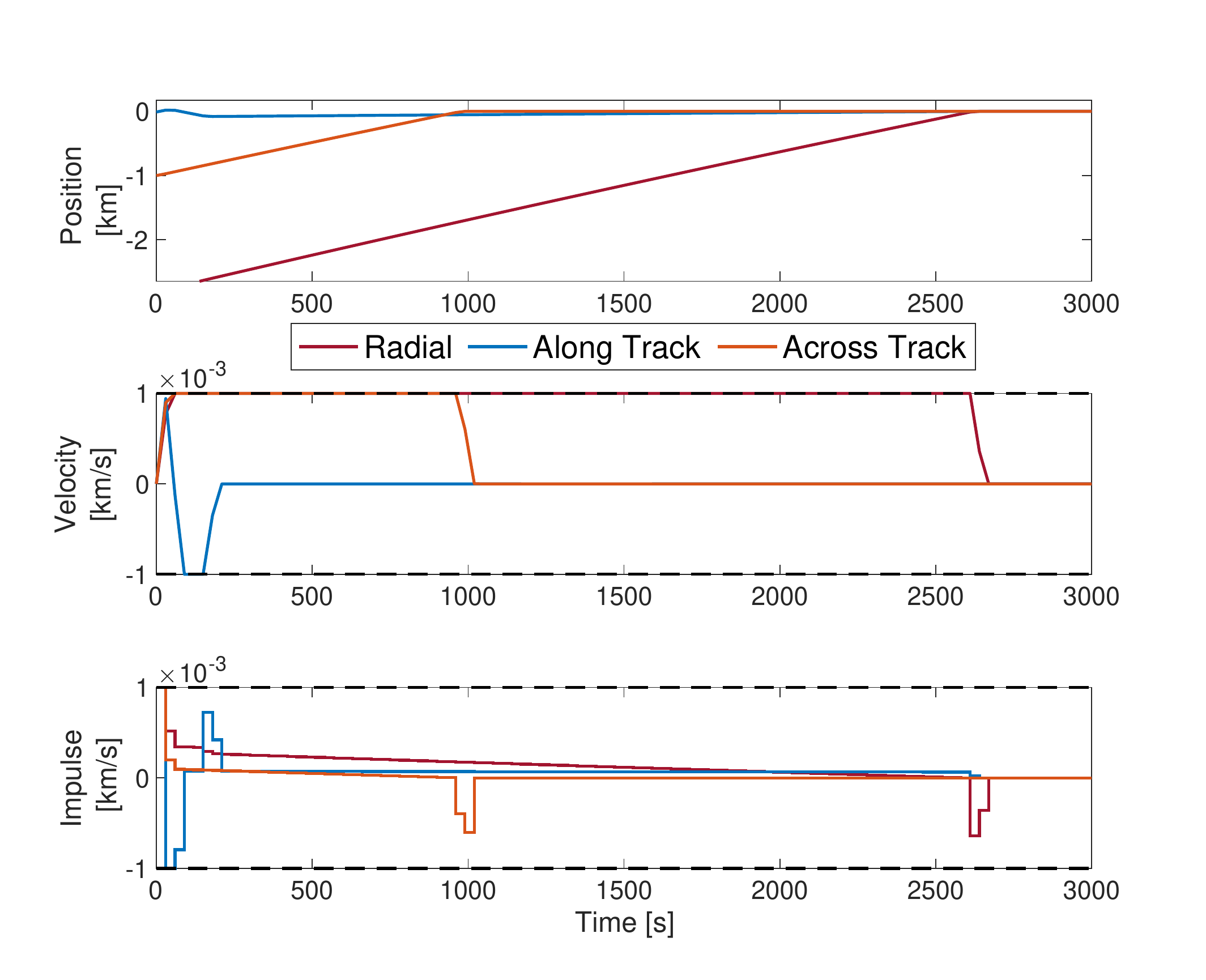}
	\caption{Closed-loop response for the spacecraft relative motion example.}
	\label{fig:hcw_cl}
\end{figure}

\textbf{Control of Spacecraft Relative Motion\cite{weiss2012model}:} These equations describe the radial, along track, and across track positions $\zeta = [x_1~x_2~x_3]^T$ and velocities $\dot{\zeta}$ of a spacecraft relative to a nominal circular orbit. The control objective is to drive the spacecraft to the origin from $\zeta_0 = -[2.8~0.01~1]~km$ and $\dot{\zeta}_0 = 0$. The system dynamics are given by the Hill-Clohessy-Wiltshire (HCW) equations,
\begin{gather}
	\ddot{x}_1 - 3\omega^2 x_1 -2 \omega \dot{x}_2 = 0,\\
	\ddot{x}_2 + 2 \omega \dot{x}_1 = 0,\\
	\ddot{x}_3 + \omega^2 x_3 = 0,
\end{gather}
where $\omega = 0.0011 ~s^{-1}$ is the mean motion of the reference orbit. The dynamics of $x = (\zeta,\dot{\zeta})$ can be compactly written as $\dot{x} = A_c x$. The control inputs are modelled as impulsive thrusts which instantaneously change the velocity of the spacecraft, see \cite{weiss2012model}, so that the discrete time model is
\begin{equation}
  x_{k+1} = A \left(\begin{bmatrix}
    \zeta_k \\ \dot{\zeta}_k
  \end{bmatrix} + \begin{bmatrix}
    0 \\ I
  \end{bmatrix} \Delta v_k\right)= A x_k + B u_k,
\end{equation}
where $u = \Delta v$ is the instantaneous change in velocity due to the impulsive thrusters, $A = e^{A_c \tau}$, and $\tau = 30~s$. The control inputs must satisfy $||u_k||_\infty \leq 1~m/s$. The spacecraft velocity is constrained to satisfy $||\dot{\zeta}_k||_\infty \leq 1~m/s$ and the tuning matrices are $Q = diag([1~1~1~0.001~0.001~0.001])$ and $R = I_{3\times 3}$. The closed-loop response of the system is shown in Figure~\ref{fig:hcw_cl}.

\textbf{Control of a Copolymerization Reactor\cite{congalidis1986modeling}:} The following normalized transfer function models the copolymerization of methyl methacrylate (MMA) and vinyl acetate (VA) in a continuous stirred tank reactor:

{\footnotesize
\begin{equation*} 
\begin{bmatrix}
\frac{0.34}{0.85s +1} & \frac{0.21}{0.42s +1} & \frac{0.25 s + 0.5}{12s^2 + 0.4s +1} & 0 & \frac{6.46(0.9s +1)}{0.07s^2 + 0.3s + 1}\\[0.8em]
\frac{-0.41}{2.41 s + 1} & \frac{0.66}{1.51s+1} & \frac{-0.3}{1.45s+1} & 0 & \frac{-3.72}{0.8s+1}\\[0.8em]
\frac{0.3}{2.54s+1} & \frac{0.49}{1.542+1} & \frac{-0.71}{1.35s+1} & \frac{-0.20}{2.71s+1} & \frac{-4.71}{0.008s^2 + 0.41s + 1}\\[0.8em]
0 & 0 & 0 & 0 & \frac{1.02}{0.07s^2 + 0.31s + 1}
\end{bmatrix}.
\end{equation*}}

\begin{figure}[htbp]
	\centering
	\includegraphics[width=0.95\columnwidth]{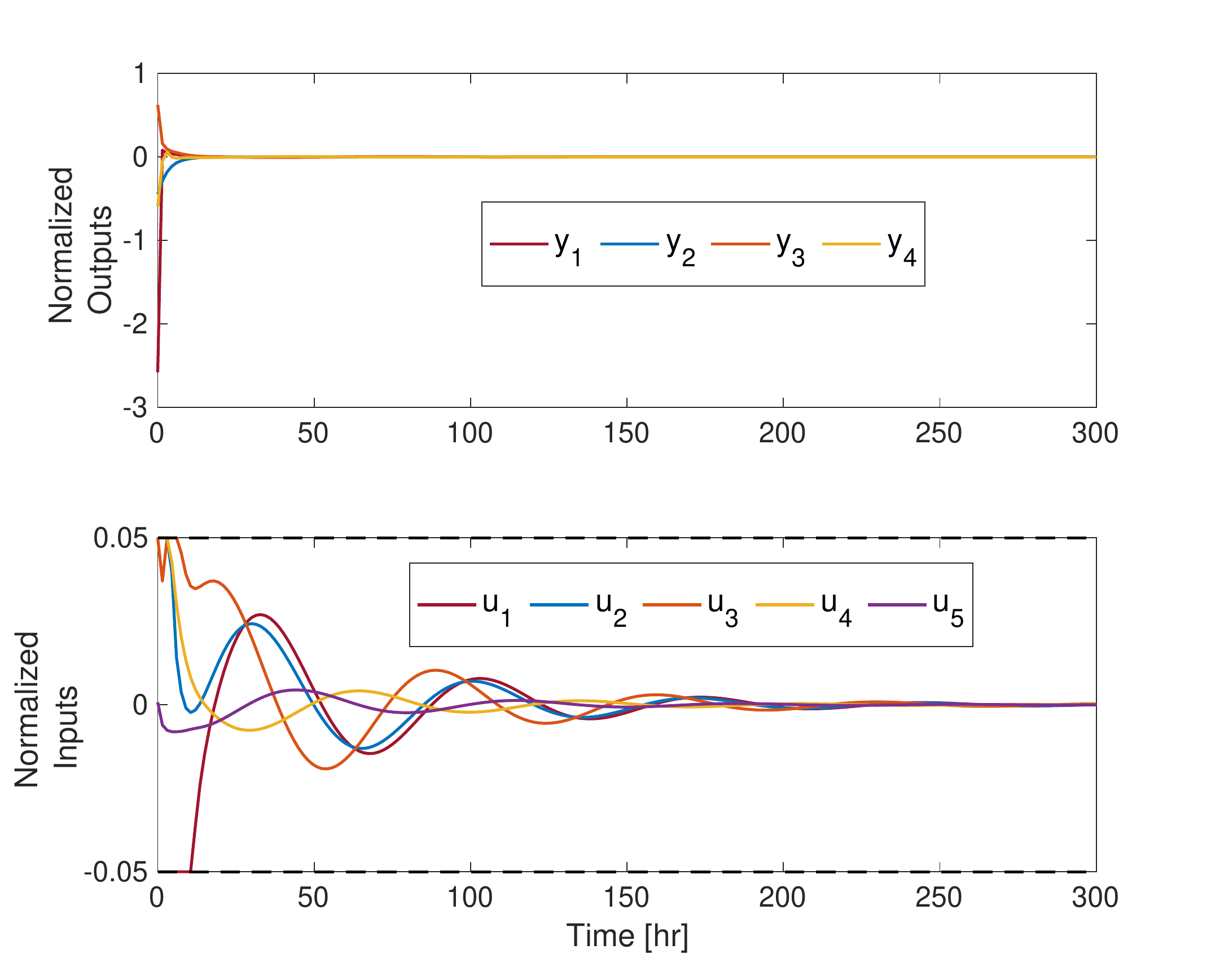}
	\caption{Closed-loop response for the copolymerization example.}
	\label{fig:copoly}
\end{figure}

The normalized inputs are flows of monomer MMA ($u_1$), monomer VA ($u_2$), initiator ($u_3$), transfer agent ($u_4$), and the reactor jacket temperature ($u_5$). The normalized outputs are the polymer production rate ($y_1$), the mole fraction of MMA in the polymer ($y_2$), the molecular weight of the polymer ($y_3$), and the reactor temperature ($y_4$). All inputs and outputs are relative to nominal operating conditions \cite{congalidis1986modeling}. The model was realized in modal form using the $\texttt{ss}$ command in MATLAB  and discretized using a zero-order hold with a normalized sampling period of $0.5$ (corresponding to three hours in physical time). The resulting model has 18 states, 5 inputs and 4 outputs. The states are initially disturbed as $\xi_{0,i} = \sin(i)$ for $i = 1,~...,~18$; the control objective is to drive the outputs to the origin. The inputs are constrained as $||u_k||_\infty \leq 0.05$, i.e., 5\% deviation from nominal. The horizon length is $N = 70$, and the weighting matrices are chosen as $Q = C^T C$, where $C$ is the output matrix from the realization process, and $R = 0.1 I_{5\times 5}$. Closed-loop traces are shown in Figure~\ref{fig:copoly}.

\begin{table}[h]
\centering
\caption{Problem data for the QPs.}
\label{tab:problem_sizes}
\resizebox{\columnwidth}{!}{
\begin{tabular}{|c|c|c|c|} \hline
 & Servo & HCW & Copoly\\\hline
 Number of States & 4 & 6 & 18 \\\hline
 Number of Controls& 1 & 3 & 5  \\\hline
 Number of timesteps & $40$ & $100$ & $200$  \\\hline
 Horizon Length & 30 &40 & 80  \\\hline
 \multicolumn{4}{|c|}{Sparse Problem} \\ \hline
 Variables  & $155$ & $369$ & $1863$  \\\hline
 Equality constraints & $124$ & $246$ & $1458$  \\\hline
 Inequality constraints& $124$ & $492$ & $810$  \\\hline
 \multicolumn{4}{|c|}{Condensed Problem} \\ \hline
 Variables  & $31$ & $123$ & $405$  \\\hline
 Inequality constraints& $124$ & $492$ & $810$  \\\hline
 Hessian condition number & $189$ & $3.19\times 10^8$ & $1.4\times 10^3$  \\ \hline
\end{tabular}}
\end{table}

\subsection{Implementation details} \label{ss:impl_details}
The default parameters used throughout are $\sigma_k = \sqrt{\epsilon_m}$, $\zeta = 10^{-14}$, $\tau_a = 10^{-4}$, $\tau_r = 0$, $\tau_{inf} = 10^{-8}$, $\alpha = 0.95$, $\beta = 0.7$, $ \eta = 10^{-8}$, where $\epsilon_m$ is machine precision; $\epsilon_m \approx 10^{-16}$ for our double precision implementation.

We have implemented two versions of FBstab \footnote{Both are available online at: \url{https://github.com/dliaomcp/fbstab-matlab.git}} in MATLAB. The first solves problems of the form \eqref{eq:QP} and is featured in in Sections~\ref{ss:infeas_tests} and \ref{ss:rt_benchmarking}. It is implemented using MATLABs built-in dense linear algebra routines and solves the Newton-step systems \eqref{eq:newton_reduced} using a Cholesky factorization. The second solves problems of the form \eqref{eq:MPCQP} and is featured in Sections~\ref{ss:scaling} and \ref{ss:rt_benchmarking}. It exploits the structure present in \eqref{eq:MPCQP} and is equipped with two different linear solvers for the Newton step systems:
\begin{enumerate}
	\item FBstab Ricatti: Uses a Ricatti-like recursion similar to the one in \cite{rao1998application} to solve \eqref{eq:newton_reduced}.
	\item FBstab MA57: Calls MA57\cite{duff2004ma57} to directly solve \eqref{eq:newton_expanded} using sparse linear algebra.
\end{enumerate}
Both implementations use the non-monotone linesearch technique of Grippo et al. \cite{grippo1986nonmonotone} to improve performance without jeopardizing the convergence properties of the algorithm. 

\begin{rmk}
In \cite{cdc_fbstab2019} we used a linear algebra framework for FBstab which solves the condensed without forming \eqref{eq:condensed_QP} explicitly. This approach uses the conjugate gradient method to solve linear systems which could be slow when \eqref{eq:condensed_QP} is ill-conditioned, e.g., in the HCW example. We found the direct methods used in paper to be faster and more robust.
\end{rmk}
\subsection{Solver scaling} \label{ss:scaling}
\begin{figure}[htbp]
	\centering
	\includegraphics[width=0.95\columnwidth]{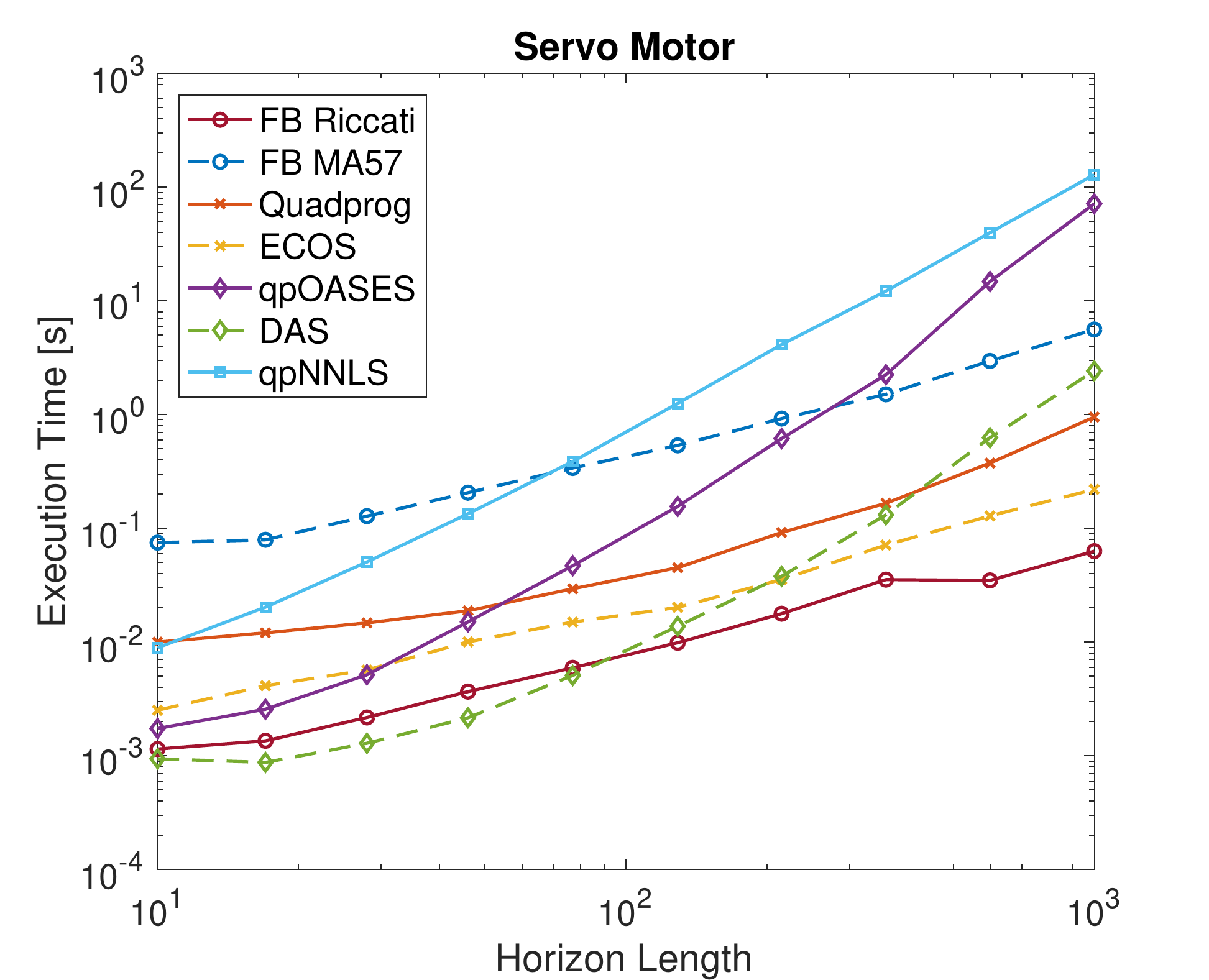}
	\caption{Solver scaling comparison for the servo motor example.}
	\label{fig:servo_scaling}
\end{figure}

\begin{figure}[htbp]
	\centering
	\includegraphics[width=0.95\columnwidth]{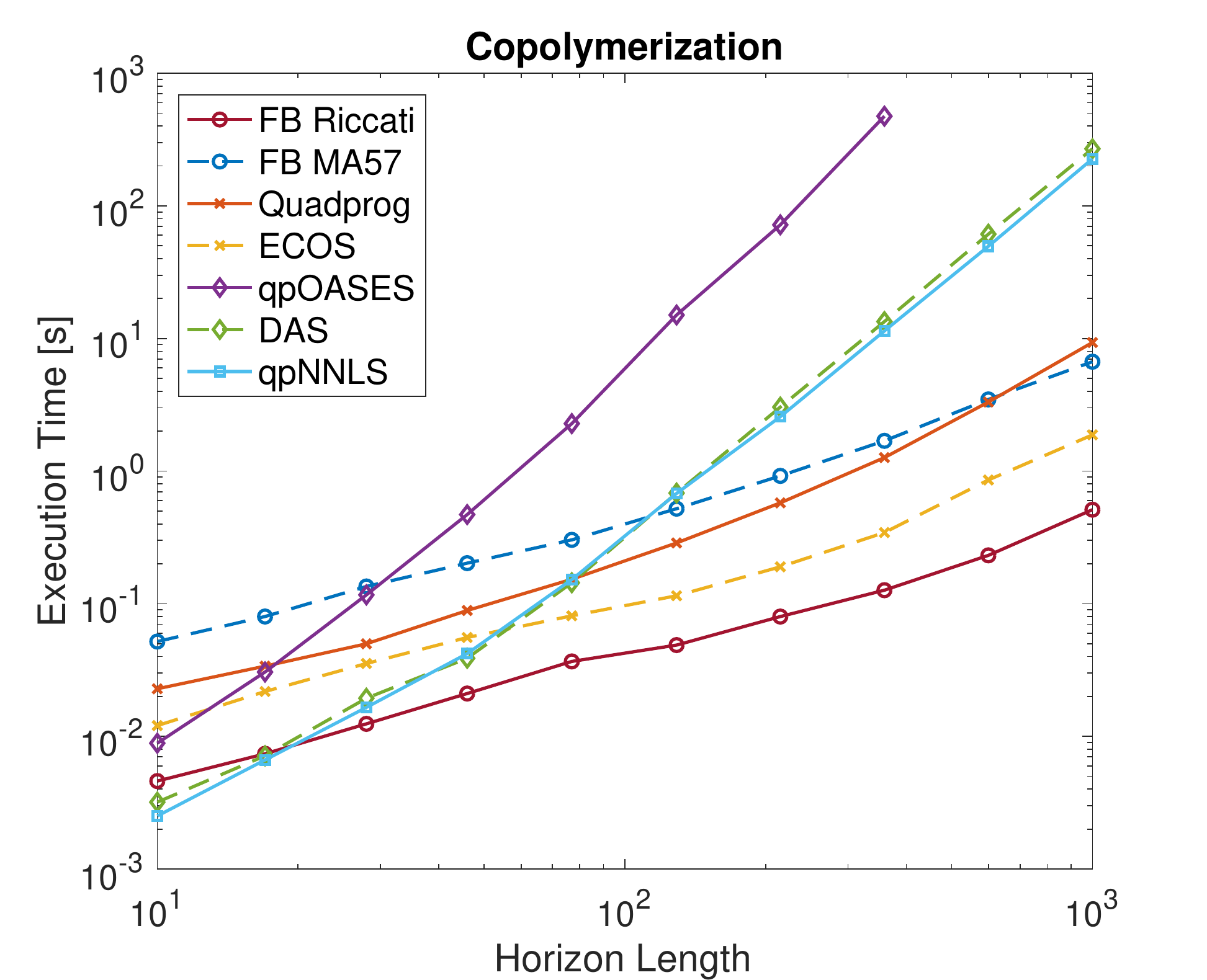}
	\caption{Solver scaling comparison for the copolymerization example.}
	\label{fig:copoly_scaling}
\end{figure}

FBstab can efficiently solve structured problems. The linear systems \eqref{eq:newton_expanded} that are solved in FBstab are highly structured similar to those in IP methods. To demonstrate this, we compared FBstab, implemented using the two different linear solvers, as described in Section~\ref{ss:impl_details}, with the external solvers: (1) \texttt{quadprog} (MATLAB 2017b) which uses sparse linear algebra and the \textit{interior-point-convex} algorithm, (2) ECOS\cite{domahidi2013ecos} (self-dual interior point) and (3) qpOASES\cite{ferreau2014qpoases} (active set). We also implemented the following in MATLAB: (4) the dual active set (DAS) method \cite{goldfarb1983numerically}, including factorization updating, (5) QPNNLS \cite{bemporad2018numerically}, (6) GPAD \cite{patrinos2014accelerated} and (7) accelerated ADMM\cite{goldstein2014fast}. The DAS method, qpNNLS, and qpOASES solve \eqref{eq:condensed_QP}, with the cost of condensing included in the analysis, while all other methods solve \eqref{eq:MPCQP} directly. The MATLAB routines were converted into \texttt{C} code using the \texttt{mex} command. We found that GPAD and ADMM were not competitive; both have been omitted from Figures~\ref{fig:servo_scaling} and \ref{fig:copoly_scaling} for clarity.

We solved the first QP, i.e., at $t = 0$, in the servo motor and copolymerization examples\footnote{The spacecraft example was omitted because the dynamics are unstable, the resulting condensed problem becomes ill-conditioned enough for large horizons to make most of the methods fail.} and measured wall clock times as the horizon was varied from $N=10$ to $N = 1000$. All methods were cold started at the origin, the experiments were performed on a 2015 Macbook Pro with a 2.8 GHz i7 processor and 16 GB of RAM running MATLAB 2017b. Recorded execution times were averaged as necessary to obtain consistent timings. Figures~\ref{fig:servo_scaling} and \ref{fig:copoly_scaling} display the results. ECOS, \texttt{quadprog} and FBstab Ricatti/MA57 scale like $\Oh(N)$ with FBstab Ricatti being the fastest method. The active set methods qpOASES, qpNNLS, and DAS are efficient for small problems but are quickly overtaken as $N$ becomes large, they scale like $\Oh(N^3)$. The DAS method is the quickest method for the servo example for short horizons before being overtaken by FBstab Ricatti. Similarly, for very short horizons, qpNNLS is the quickest method for the copolymerization example before being overtaken by FBstab. As expected, active set methods are very effective for small problems but are quickly overtaken by both interior point methods and FBstab. FBstab is shown to be faster than several interior point methods at all horizon lengths. Overall, FBstab scales well as is competitive with and often superior to several established solvers. 

\subsection{Benchmarking on real-time hardware} \label{ss:rt_benchmarking}
To investigate the performance of FBstab on embedded hardware we performed some benchmarking on a Speedgoat Baseline Real-time Target Machine (SGRTM). The SGRTM (2.0 GHz Celeron CPU, 4 GB RAM) is a rapid prototyping platform which runs a real-time operating system (RTOS) and is representative of an embedded computing environment. Using the RTOS allows us to obtain deterministic execution time measurements. Moreover, the SGRTM does not have any linear algebra libraries available, e.g., BLAS or LAPACK, so solving ill-conditioned problems is more difficult.

We implemented FBstab using structure exploiting linear algebra (FBstab MPC) and using dense linear algebra (FBstab Dense) as described in Section~\ref{ss:impl_details}. For comparison we implemented the following in MATLAB: (1) FBRS (Fischer-Burmeister Regularized and Smoothed)\cite{liaomcphersonFBRS}, (2) DAS (Dual Active Set) \cite{goldfarb1983numerically}, (3) QPNNLS (robust nonnegative least squares\cite{bemporad2018numerically}), (4) PDIP (primal-dual IP using Mehrotra's predictor-corrector \cite[Algorithm 14.3]{nocedal2006numerical}), (5) GPAD (accelerated dual gradient projection)\cite{patrinos2014accelerated}, and (6) accelerated ADMM \cite{goldstein2014fast}. Since \texttt{quadprog} cannot codegenerate, ECOS does not have a Simulink interface, and qpOASES was outperformed by DAS during the scaling trials all three were omitted from testing; GPAD and ADMM were not competitive and were omitted as well. The methods were converted into C code using Simulink Real-time (2017b). A method is deemed to have failed if it is stalls and is unable to solve any of the QPs in the sequence to the desired precision ($10^{-4})$. Note that we have implemented $LDL^T$ and $QR$ factorization updating for qpNNLS and DAS to ensure a competitive comparison.

The results are shown in Table~\ref{tab:sg_warm}. When warmstarting is enabled, FBstab MPC is the fastest method in the worst case for all three examples and FBstab Dense is competitive with the IP and AS methods, especially on the larger copolymerization example. In terms of average execution times\footnote{Average execution times are an indicator of power draw. This is an important metric in Aerospace applications where reduced power consumptions leads to e.g., extended range for drones.} FBstab MPC and FBstab Dense are dominant. When warmstarting is disabled, the PDIP and active set methods become more competitive, however FBstab MPC is still more efficient.

Overall, when the cost of condensing is considered, see Remark~\ref{rmk:COC}, FBstab MPC is shown to outperform the other methods tested in terms of both maximum and average execution time. FBstab derives significant benefit from warmstarting, this is especially noticeable for the HCW example, and is significantly faster than FBstab Dense, showcasing the importance of specialized linear algebra routines. Further, both DAS and FBRS fail on the ill-conditioned HCW example while their regularized versions, qpNNLS and FBstab respectively, succeed, demonstrating the expected improved robustness due to proximal regularization. FBstab is often faster than FBRS, demonstrating that the addition of proximal regularization makes the methods more robust without a significant reduction in speed.

\begin{rmk} \label{rmk:COC}
The cost of condensing, i.e., of converting \eqref{eq:MPCQP} to \eqref{eq:condensed_QP} is included in the results reported in Table~\ref{tab:sg_warm}. This simulates solving e.g., trajectory tracking problems or real-time iteration \cite{diehl2005real} subproblems. These computations can sometimes be moved offline, in this situation the normalized cost of condensing listed in Table~\ref{tab:sg_warm} should be subtracted from each row of the last five columns. In this scenario, the DAS method is the best method for the servo motor example.
\end{rmk}

\begin{table}[h]
\caption{Summary of normalized Speedgoat benchmarking reporting the maximum and average QP solutions times for each sequence. Warm and cold starting are indicated by W and C respectively.}
\label{tab:sg_warm}
\resizebox{\columnwidth}{!}{
\begin{tabular}{|c|c|c|c|c|c|c|} \hline
 & FBstab & FBstab & FBRS & PDIP & NNLS & DAS\\
 & MPC & Dense & &  & & \\ \hline
\multicolumn{7}{|c|}{Servo Motor, Normalization $= 4.5~ms$} \\ 
\multicolumn{7}{|c|}{Normalized Cost of Condensing $= 1.1$} \\ \hline
MAXW & 1.00 & 2.9 & 2.4 & 3.2 & 3.3 & 1.7 \\ \hline
AVEW & 0.2 & 1.4 & 1.4 & 2.1 & 1.8 & 1.3  \\ \hline
MAXC & 1.5 & 5.0 & 5.0 & 3.0 & 2.9 & 2.1  \\ \hline
AVEC & 0.2 & 1.6 & 1.6 & 2.1 & 1.8 & 1.3  \\ \hline
\multicolumn{7}{|c|}{Spacecraft, Normalization $= 63.9~ms$} \\ 
\multicolumn{7}{|c|}{Normalized Cost of Condensing $= 1.5$} \\ \hline
MAXW & 1.00 & 14.8 & F & 11.7 & 15.6 & F \\ \hline
AVEW & 0.1 & 3.3 & F & 8.3 & 6.5 & F  \\ \hline
MAXC & 3.4 & 73.2 & F & 29.3 & 7.6 & F  \\ \hline
AVEC & 2.2 & 62.8 & F & 25.5 & 3.7 & F  \\ \hline
\multicolumn{7}{|c|}{Copolymerization, Normalization $= 97.1~ms$} \\ 
\multicolumn{7}{|c|}{Normalized Cost of Condensing $= 76.3$} \\ \hline
MAXW & 1.00 & 96.6 & 102.9 & 238.7 & 94.4 & 149.2 \\ \hline
AVEW & 0.4 & 82.5 & 82.6 & 204.2 & 85.8 & 96.2 \\ \hline
MAXC & 1.5 & 113.5 & 112.9 & 238.2 & 88.4 & 293.3  \\ \hline
AVEC & 0.3 & 83.0 & 82.9 & 205.3 & 85.9 & 101.7 \\ \hline
\end{tabular}}
\end{table}

\subsection{Degenerate and Infeasible Problems} \label{ss:infeas_tests}
Consider the following parameterized QP:
\begin{subequations}\label{eq:infeas_testingQP}
\begin{gather} 
\underset{x_1,x_2}{\mathrm{min.}} \quad \frac12 x_1^2 + x_1 + c x_2\\
\mathrm{s.t} \quad a_1 x_1 + a_2 x_2 \leq 0,\\
1 \leq x_1 \leq 3,\\
1 \leq x_2 \leq b,
\end{gather}
\end{subequations}
by varying $a_1,a_2,b$ and $c$ we can create degenerate or infeasible test problems.

First, we consider degeneracy. Setting $a_1 = a_2 = c = 0$ and $b = 3$ we obtain a degenerate QP with the primal solution set $\Gamma_p = \{1\}\times [1,3]$. FBstab signals optimality after 2 proximal iterations and 5 Newton iterations and returns $x^* = (1.00,1.00), v^* = (0.00,0.00,0.00,2.00,0.00)$ with the norm of the residual $\epsilon = 8.97\times 10^{-12}$. Second, we consider a primal infeasible QP by setting $a_1 = a_2 = 0, c= -1$ and $b = 3$. FBstab signals primal infeasibility after 1 proximal iteration and 7 Newton iterations and returns $\delta x^* = (-0.18,0.36)\cdot10^{-5}, \delta v^* = (4.47,0,0,4.47,4.47)\cdot10^7$. Finally, we consider a dual infeasible QP by setting $a_1 = a_2 = 0, c = -1$ and $b = \infty$.  This leads to a QP for which $x = (0,1)$ is a direction of unbounded descent. FBstab signals dual infeasibility after 3 proximal iterations and 8 Newton iterations and returns $\delta x = (0,671)\cdot 10^5$ and $\delta v = 0$. FBstab was initialized at the origin for all results reported in this section.

\section{Experimental Results}

\begin{figure}[htbp]
  \centering
  \includegraphics[width=0.95\columnwidth]{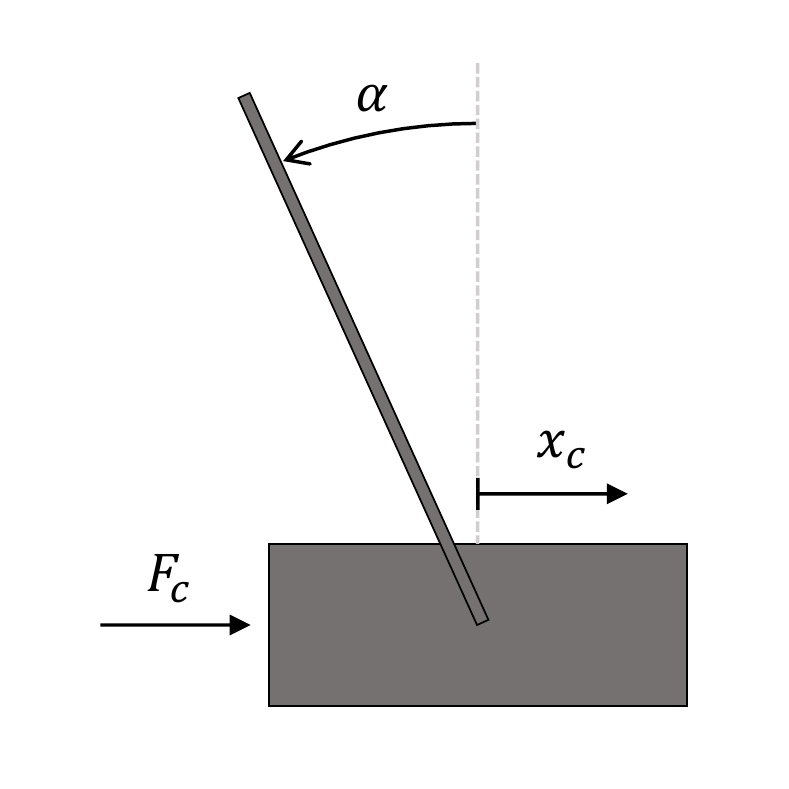}
  \caption{Diagram of the inverted pendulum on a cart system.}
  \label{fig:pen_cart_diagram}
\end{figure}

In this section we showcase FBstab's utility with an experimental demonstration. Consider a pendulum mounted on a cart as shown in Figure~\ref{fig:pen_cart_diagram}. The control objective is to drive the cart position, $x_c$, to a target value while balancing the pendulum, i.e., keeping $\alpha \approx 0$. The nonlinear equations of motion of the cart-pendulum system are given in \eqref{eq:cart_eom}, a motor is used to accelerate the cart. The corresponding parameters are given in Table~\ref{tab:pen_param}; these values correspond to a QUANSER \textit{Linear Servo Base Unit with Inverted Pendulum} device. The states and control inputs are
\begin{equation}
  x = (x_c,\alpha,\dot{x}_c,\dot{\alpha}), ~ u = V_m,
\end{equation}
i.e., the cart position, pendulum angle, their velocities, and the motor input voltage.

\begin{subequations} \label{eq:cart_eom}
\begin{gather}
a\ddot{x}_c - b\cos(\alpha) \ddot{\alpha} + b \sin(\alpha) \dot{\alpha}^2 + c\dot{x} = F_c, \\
-b \cos(\alpha) \ddot{x}_c + d \ddot{\alpha} - e\sin(\alpha) + f \dot{\alpha} = 0,\\
F_c = -k_1 \dot{x}_c + k_2 V_m
\end{gather}
\end{subequations}

\begin{table}[h]
\caption{Inverted Pendulum Parameters}
\label{tab:pen_param}
\begin{tabular}{|c|c|c|c|} \hline
Parameter & Value & Parameter & Value \\ \hline
$a $& 1.3031 & $b$ & 0.0759 \\ \hline
$c$ & 5.400 & $d$ & 0.0330 \\ \hline
$e$ & 0.7450 & $f$ & 0.0024 \\ \hline
$k_1$ & 7.7443 & $k_2$ & 1.7265 \\ \hline
\end{tabular}
\end{table}

We implemented a linear MPC controller of the form \eqref{eq:MPCQP} with $N = 10$. The prediction model is obtained by linearizing \eqref{eq:cart_eom} about the origin then discretizing the resulting linear continuous time model using a zero order told at a sampling time of $0.01 s$. The weighting matrices are $Q_i  = Q = diag(35,0.01,0.1,0.1),~i \in \intrng{0}{N-1} $, $R_i = R = 0.02,~i \in \intrng{0}{N}$, $S_i = 0,~i\in \intrng{0}{N}$ and $Q_N = P$, where $P$ is the solution of the discrete time algebraic Ricatti equation corresponding to the prediction model. We also impose the following constraints on the input voltage and the pendulum angle $|\alpha| \leq 4^\circ,~|u| \leq 10~V$.

We implemented the controller in MATLAB/SIMULINK using the QUANSER QUARC real-time software to interface with the experimental setup; the system sampling rate is $100$ Hz. We solved \eqref{eq:condensed_QP} using the FBstab Dense, the pendulum angle constraint was softened using $L_2$ penalties to ensure feasibility despite model mismatch. The resulting QP has $22$ decision variables and $55$ inequality constraints, each QP was warmstarted with the solution from the previous timestep.

\begin{figure}[htbp]
  \centering
  \includegraphics[width=0.95\columnwidth]{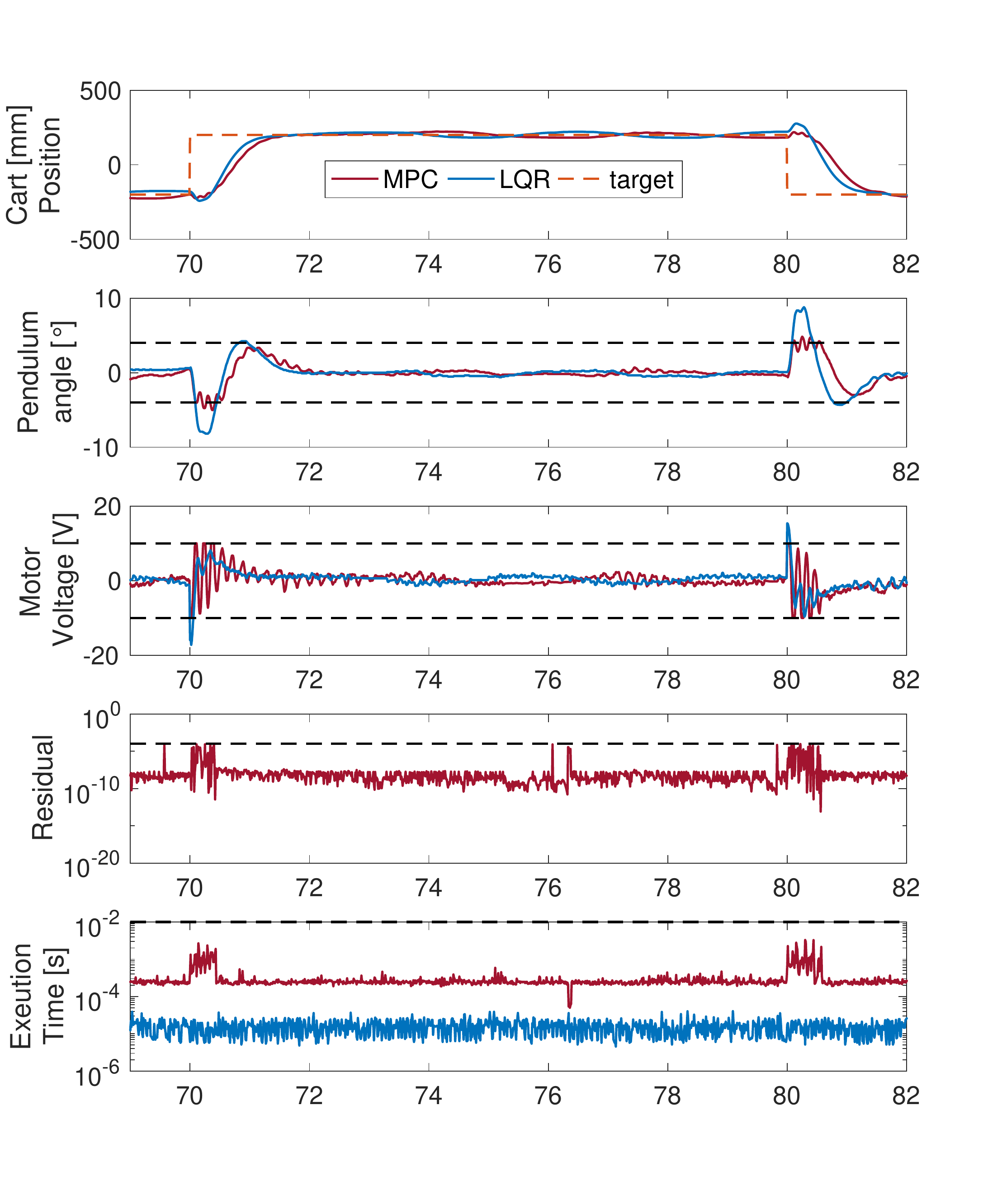}
  \caption{Experimental stabilization of an inverted pendulum on a cart using MPC and a LQR. The oscillations are due to unmodelled flexibility in the track the cart runs along and of the table the system is mounted on.}
  \label{fig:exp_cart}
\end{figure}

The experimental results are shown in Figure~\ref{fig:exp_cart}. The MPC controller is able to track the desired cart position trajectory while balancing the pendulum and respecting input constraints. The pendulum angle constraint is enforced with some minor violations caused by unmodelled dynamics\footnote{Specifically, vibration of the table the device is mounted on and flexibility of the track the cart rides along.}. FBstab successfully maintains the QP residual below the prescribed precision of $10^{-4}$. The maximum execution time (measured using the QUARC software) was $4.3~ms$, safely below the sampling period of $10~ms$.

\section{Conclusions}
This paper presents FBstab, a proximally stabilized Fischer-Burmeister method for convex quadratic programming. FBstab is attractive for real-time optimization because it is easy to code, numerically robust, easy to warmstart, can exploit problem structure, and converges or detects infeasibility under only the assumption that the Hessian of the quadratic program is convex. An open source MATLAB implementation of FBstab is available online. Future work includes exploring the application of stabilized semismooth Newton-type methods to nonlinear problems and preparation of an open source \texttt{C++} implementation of FBstab.

\section*{Acknowledgments}
We would like to thank Dr. Chris Petersen of the U.S. Air Force Research Laboratories for encouraging us to develop methods for optimization problems with degenerate solutions and Brian Ha for helping to set up the inverted pendulum system.

\bibliography{fbstab}

\end{document}